\documentclass[journal]{IEEEtran}

\usepackage{myStyleIEEE}
\usepackage[margin=0.52in]{geometry}

\DeclareSIUnit{\sieuro}{\mbox{\euro}}
\makenomenclature

\newcolumntype{L}[1]{>{\raggedright\let\newline\\\arraybackslash\hspace{0pt}}m{#1}}
\newcolumntype{C}[1]{>{\centering\let\newline\\\arraybackslash\hspace{0pt}}m{#1}}
\newcolumntype{R}[1]{>{\raggedleft\let\newline\\\arraybackslash\hspace{0pt}}m{#1}}

\IEEEoverridecommandlockouts

\begin{document}

\title{Economic Valuation and Pricing of Inertia in Inverter-Dominated Power Systems}

\renewcommand{\theenumi}{\alph{enumi}}

\newcommand{\uros}[1]{\textcolor{magenta}{$\xrightarrow[]{\mathrm{UM}}$ ``#1''}}
\newcommand{\vaggelis}[1]{\textcolor{blue}{$\xrightarrow[]{\mathrm{EV}}$ ``#1''}}
\newcommand{\SD}[1]{\textcolor{red}{ #1}}
\newcommand{\UM}[1]{\textcolor{magenta}{ #1}}
\newcommand{\matthieu}[1]{\textcolor{pPurple}{$\xrightarrow[]{\mathrm{MP}}$ ``#1''}}

\author{Matthieu~Paturet, 
        Uros~Markovic,~\IEEEmembership{Member,~IEEE,}
        Stefanos Delikaraoglou,~\IEEEmembership{Member,~IEEE,}
        Evangelos Vrettos,~\IEEEmembership{Member,~IEEE,}
        Petros Aristidou,~\IEEEmembership{Senior~Member,~IEEE,}
        and~Gabriela~Hug,~\IEEEmembership{Senior~Member,~IEEE}
        \vspace{-0.5cm}}

\maketitle
\IEEEpeerreviewmaketitle

\begin{abstract}
This paper studies the procurement and pricing of inertial response using a frequency-constrained unit commitment formulation, which co-optimizes the provision of energy and inertia services while accounting for their complementary properties. The proposed approach builds on a two-step process that allows to differentiate between the units being online for energy purposes and the ones committed additionally solely for inertia provision. Subsequently, three novel pricing and payment schemes that reimburse inertia providers in a transparent and fair manner according to their individual participation are proposed. The analysis considers both synchronous and converter-based generators and provides insights regarding the impact of each pricing scheme on total system cost, as well as on the investment signals for technologies capable of offering affordable inertial response. The results show that all three methods will have a beneficial impact on frequency stability and aid the operator in ensuring system reliability, whereas the policy implications to different inertia providing units will vary between the payment schemes. 
\end{abstract}


\vspace{-0.35cm}
\section{Introduction} \label{sec:Intro}

\lettrine[lines=2]{C}{urrent} strategies for power system operation and the existing market architecture were designed in a context where conventional generators (e.g., hydro, coal or gas units) were dominating. The spinning turbines of these generators are synchronized to the grid frequency and their large rotating masses naturally provide the necessary physical inertia to contain frequency deviations caused by sudden power imbalances. Nowadays however, increasing shares of Renewable Energy Sources (RES) with null marginal cost are displacing the more expensive Synchronous Generators (SGs). Considering that RES units are typically interfaced to the grid via power electronics and therefore do not provide physical inertia to the network, the ability of the power system to dampen frequency changes is diminished. In case of insufficient inertial response, large disturbances will result in fast frequency dynamics and high Rate-of-Change-of-Frequency (RoCoF), which may even lead to disconnection of generators or loads \cite{Milano2018}.

In view of this new operational paradigm, inertia cannot be viewed anymore just as byproduct of the power provision from synchronous generators. In addition, existing frequency control mechanisms such as droop response by frequency containment reserves are not able to instantly reduce frequency changes following a disturbance~\cite{Hartmann2019}. To this end, it becomes imperative to design new market mechanisms and tailored market products that are able to adequately reward the provision of Fast Frequency Response (FFR) services from existing sources, i.e., synchronous generators, and provide correct price signals, that will incentivize new investments in units that can contribute to the inertial response of the power system such as fast frequency-responsive loads (e.g., heat pumps, batteries) or synthetic inertial response from wind generation. In the meantime, the new market mechanisms should allow the system operator to express its willingness to pay for these services dynamically, depending on system conditions and the scarcity of FFR services.


Motivated by the above considerations, the system operator in Great Britain (National Grid) has designed new frequency regulation products to bring more fast-acting assets into the system \cite{Nationalgrid}. In the same vein, a new FFR reserve will be rolled out in the Nordic synchronous area \cite{FFRnordic} in spring of 2020. These are distinct products tailored to enhance the fast frequency response of the system. However, being traded in separate markets, i.e., unbundled from energy, may fail to capture the complementarity between inertia and energy services and consequently increase the total system cost.

Following a different approach, Eirgrid in Ireland has introduced a Synchronous Inertial Response (SIR) constraint to ensure that the available inertia does not fall below a static limit of \SI{23000}{\mega\watt\second} \cite{EirGridReport3}. However, as this mechanism remunerates service providers based on an annual tariff, it does not reflect the instantaneous value of inertia for the system. On the other hand, the system operator in Texas (ERCOT) performs a dynamic dimensioning of Responsive Reserve Services (RRS) that are used to maintain the reliability of the grid during low-inertia situations \cite{ERCOTwind}. Since energy and RRS are co-optimized, the resulting prices account for the fact that these products compete for the same generation capacity. Nonetheless, RRS bundle various services (load resources with under-frequency relays, governor response and frequency restoration provided by generators after a frequency event) that have different contributions to frequency response and are not necessarily reflective of inertia provision. This makes the accurate definition of RRS requirements even more complex, specifically in terms of inertia, and it dilutes the price signals to the assets that truly contribute to inertial response.

Previous works in \cite{xu2016economic} and \cite{badesa2017economic} have quantified the value of inertia as the decrease in total cost of a frequency-constrained Unit Commitment (UC) model when including cost-free synthetic inertia. These studies showcase the benefits in terms of operational cost from enlarging the pool of inertia providers, but they disregard the costs of converter-interfaced generators and do not discuss any pricing issues. This is a non-trivial problem as the discrete nature of inertia, depending on the on/off status of the SGs, introduces pricing complexities when integrating frequency constraints. In order to circumvent the underlying issues, several studies use a linearized version of these constraints. A marginal pricing scheme is proposed in \cite{ela2014market1} based on linear constraints fitted to simulation data. Authors in \cite{greve2018system} focus on pricing frequency response services of varying quality (i.e., with different response times) using a two-sided Vickrey-Clarke-Groves auction, where the preferences of the system operator for each service are reflected in different utility functions. Finally, the work in \cite{badesa2019pricing} presents a marginal pricing scheme for frequency services based on a Mixed-Integer Second-Order Cone Program (MISOCP) formulation of the frequency-security constraints. However, while the resulting prices provide reasonable incentives to the inertia providers, they do not guarantee cost recovery of all synchronous generators. Moreover, the frequency response from different unit types is characterized solely by the response time, which fails to capture the fundamental differences in inertia provision from conventional and converter-based generators.

In this work, we study the procurement and pricing of inertial response using a frequency-constrained UC problem, which co-optimizes the provision of energy and inertia services accounting endogenously for their complementary properties. By incorporating the analytic expressions for frequency constraints in our UC model, we perform a dynamic dimensioning of FFR reserve requirements instead of using exogenously defined static requirements. The proposed approach builds on a two-step process, akin to the reliability UC mechanism that is already implemented in various markets in the United States \cite{gribik2007market}. Such arrangement allows to differentiate between the units being online for energy purposes and the ones committed additionally solely for inertia provision. We contribute to the ongoing discussion on FFR market design by proposing three pricing and payment schemes that reimburse inertia providers in a transparent and fair manner according to their individual participation. We analyze how these pricing schemes can lower system cost by attracting new investments in technologies capable of offering inertia at lower cost, such as grid-supporting converter control schemes that employ droop and Virtual Inertia (VI) algorithms providing synthetic inertia and Fast Frequency Control (FFC) services at a more competitive cost than SGs \cite{ofir2018droop}.

The rest of the paper is structured as follows. In Section~\ref{sec:FreqModel}, we derive the post-contingency frequency dynamics and the corresponding analytic expressions are incorporated into the UC model. The economic valuation of inertia provision and the notion of inertia substitution are explained in Section~\ref{sec:EconVal}. Section~\ref{sec:InertiaPricing} presents the three proposed inertia pricing schemes and the corresponding payments to participating units. Section~\ref{sec:CaseStudy} analyzes the performance of each pricing method on an illustrative system and on a larger test case with and without the inclusion of virtual inertia. Finally, Section~\ref{sec:Conclusion} draws the main conclusions and discusses the outlook of the study. 

\vspace{-0.15cm}
\section{Frequency-Constrained Unit Commitment} \label{sec:FreqModel}

\subsection{Frequency Dynamics in a Low-Inertia System} \label{subsec:FreqDynam}

We employ a frequency response model of a generic low-inertia system previously proposed and described in \cite{StochUCinLowInertiaGrids,UrosLQR} that comprises both traditional and converter-interfaced generators. The generator dynamics are modeled by the swing equation and governor control, whereas the converter control schemes encompass both droop and virtual inertia-based approaches, being the two of the currently most prevalent grid-supporting control techniques in the literature \cite{Rocabert2012,Tamrakar2017,UrosGM}. By employing a set of trivial mathematical operations and assumptions on the parameter selection, the time-domain evolution of the frequency deviation after a step disturbance $\Delta P\in\R$ can be computed analytically, which subsequently yields the expressions for frequency metrics of interest such as frequency nadir ($\Delta{f}_\mathrm{max}\in\R$), RoCoF ($\dot{f}_\mathrm{max}\in\R$) and steady-state frequency deviation ($\Delta f_\mathrm{ss}\in\R$):
\begin{subequations} \label{eq:puConst}
\begin{align}
    \Delta f_\mathrm{max} &= - \frac{\Delta P}{D+R_g} \left( 1 + \sqrt{\dfrac{T(R_g-F_g)}{M}} e^{-\zeta\omega_n t_m} \right), \label{eq:nadir}\\
    \dot{f}_\mathrm{max} &= \dot{f}(t_0^+) = -\frac{\Delta P}{M}, \;\;\Delta f_\mathrm{ss} = -\frac{\Delta P}{D + R_g}. \label{eq:rocof}
\end{align}
\end{subequations}
In \eqref{eq:nadir}, $\zeta\in\R_{>0}$ is the damping ratio, $\omega_n\in\R_{>0}$ is the natural frequency and $t_m\in\R_{\geq0}$ is the time instance of frequency nadir. More details on the mathematical formulation, variable definitions and model validation can be found in \cite{UrosLQR}. 

It can be observed that the relevant frequency metrics are directly dependent on the aggregate system parameters, namely total inertia $M\in\R_{>0}$, damping $D\in\R_{>0}$, droop gain $R_g\in\R_{\geq0}$ and fraction of total power generated by the high-pressure turbines of synchronous machines $F_g\in\R_{\geq0}$, and thus they could be influenced through UC decisions. In particular, RoCoF and steady-state deviation can be explicitly controlled via $\dot{f}_{\mathrm{max}}\sim M^{-1}$ and $\Delta f_\mathrm{ss}\sim (D+R_g)^{-1}$, while nadir can be modeled using a highly nonlinear function $\Delta f_{\mathrm{max}}\left(M,D,R_g,F_g\right)$.

\vspace{-0.35cm}
\subsection{RoCoF Constraint for Inertia Valuation and Pricing}

With special interest in inertia, the most relevant constraint of the three presented in \eqref{eq:puConst} is the limit on maximum instantaneous RoCoF in \eqref{eq:rocof}, with the aggregate inertia constant $M$ being the only decision variable. The quasi steady-state frequency deviation is affected solely by the damping and generator droop constants, whereas the frequency nadir depends on all aforementioned variables. In our previous work in \cite{StochUCinLowInertiaGrids}, we have observed that violating the RoCoF threshold tends to be drastically more common than reaching the critical nadir limit, thus suggesting that the RoCoF constraint is generally the binding one and therefore the most relevant when looking at the valuation and pricing of inertia. Moreover, several inertia market reports from different Transmission System Operators (TSOs) have focused solely on RoCoF, claiming that its impact on inertia levels is more important than the one of nadir \cite{EnergiforskReport, AEMCreport}. Indeed, RoCoF is an instantaneous metric related to islanding and protection scheme issues, whereas nadir affects mostly load shedding and can be regulated in real-time using FFC provision from converter-interfaced generation. Given the above reasons, this work focuses on pricing inertia to meet the RoCoF constraint requirements.


\vspace{-0.35cm}
\subsection{Unit Commitment Formulation} \label{sec:UC_formulation}

Here, we first provide the mathematical formulation of the frequency-constrained UC model comprising only synchronous generators, and subsequently extend it by including the VI units.

\subsubsection{UC without Virtual Inertia}
Let us denote by $\mathcal{T}$ the index set of different time periods and by $\mathcal{I}$, $\mathcal{J}$ and $\mathcal{N}$ the index sets of conventional generators, wind farms and nodes in the network, respectively. The mathematical formulation of the UC model without VI provision is given by
\begin{subequations}
\label{eq:UCopt1}
\begin{alignat}{3}
    & \min_\Phi \quad && \sum_{t\in \mathcal{T}} \sum_{i\in \mathcal{I}} \big(C^{\mathrm{SU}}_{i}y_{it} + C_{i} p_{it}\big) \label{eq:UC_obj}\\
    & \;\mathrm{s.t.} \quad && \nonumber \\
    & && \sum_{i \in \mathcal{I}} p_{it} +  \sum_{j \in \mathcal{J}} W_{jt} - L_{nt} = 0, \; \forall n,t,  \label{eq:pow_bal_da}\\
    & && u_{it} - u_{i(t-1)} \leq u_{i\tau^{\mathrm{1}}},  \;   \forall i,t\geq1, \label{eq:min_on}\\ 
    & && u_{i(t-1)} - u_{it} \leq 1 - u_{i\tau^{\mathrm{0}}}, \; \forall i,t\geq1,  \label{eq:min_off}\\
	& && y_{it} \geq u_{it} - u_{i(t-1)}, \; \forall i,t\geq1 \label{eq:startup}\\
	& && \ushort{P}_i u_{it} \le p_{it} \leq \widebar{P}_i u_{it}, \; \forall i,t, \label{eq:pmax}\\
	& && k_{i t} = \frac{\widebar{P}_i}{\sum_{i\in \mathcal{I}}  \widebar{P}_i } u_{it},  \; \forall i,t, \label{eq:iner_K}\\
	& && M_{t} = \sum_{i\in \mathcal{I}} 2 H_i k_{i t},  \; \forall i,t, \label{eq:iner_M}\\
	& && \frac{\dot{f}_\mathrm{lim}}{f_0} M_{t} \geq \lvert\Delta P_{t}\rvert,   \;  \forall t,  \label{eq:rocof_cstr} \\
    & && p_{it},k_{it} \ge 0,  \forall i, t ; \delta_{n t} \ge 0, \forall n,t ; M_t \ge 0, \forall t; \nonumber \\
    & && u_{it}, y_{it} \in \{0,1\}, \forall i,t, \label{eq:vardef} 
\end{alignat}
\end{subequations}
where $\Phi = \{p_{it},u_{it},y_{it},k_{it}, \forall i,t; \delta_{nt}, \forall n,t; M_{t}, \forall t \}$ is the set of optimization variables. 

The objective function \eqref{eq:UC_obj} to be minimized is the total system cost that comprises the fuel and the start-up costs, denoted by $C_{i}$ and $C^{\mathrm{SU}}_{i}$, respectively. Equality constraint \eqref{eq:pow_bal_da} enforces the nodal power balance, with $p_{it}$ being the power output of SG unit $i\in\mathcal{I}$, $W_{jt}$ referring to the power output of wind farm $j\in\mathcal{J}$, and $L_{nt}$ representing the demand at node $n\in\mathcal{N}$. Constraints \eqref{eq:min_on}-\eqref{eq:min_off} model the minimum online and offline time of conventional units based on commitment variable $u_{it}\in\{0,1\}$, while parameters $\tau^{\mathrm{1}}_g$ and $\tau^{\mathrm{0}}$ are defined as $\tau^{\mathrm{1}} = \min \{t+\mathcal{T}^{\mathrm{1}}_g-1,\ T\}$ and $\tau^{\mathrm{0}} = \min \{t+\mathcal{T}^{\mathrm{0}}_g-1,\ T\} $, where $\mathcal{T}^{\mathrm{1}}_g$ and $\mathcal{T}^{\mathrm{0}}_g$ denote the duration of time for which unit $i$ should remain online and offline, respectively. Constraint \eqref{eq:startup} models the start-up of conventional units using the binary variables $y_{it}\in\{0,1\}$. The scheduled energy production is bounded by the minimum and maximum generation limits in \eqref{eq:pmax}. The set of constraints \eqref{eq:iner_K}-\eqref{eq:rocof_cstr} captures the frequency dynamics of the system and imposes RoCoF thresholds, with $f_0=\SI{50}{\hertz}$ being the nominal frequency. In particular, equality \eqref{eq:iner_K} defines $k_{i t}$ as the ratio of the installed power capacity of generator $i$ and the total system capacity multiplied by the binary unit commitment variable $u_{it}$. Expression \eqref{eq:iner_M} introduces the average system variable for inertia, while \eqref{eq:rocof_cstr} enforces the RoCoF limit dependent on the predefined parameter $\dot{f}_\mathrm{lim}$ prescribed by the system operator \cite{entsoe}. Finally, \eqref{eq:vardef} declares the decision variables.

\subsubsection{UC with Virtual Inertia}

Let us denote by $\mathcal{V}\subseteq\mathcal{J}$ the subset of wind farms equipped with batteries for provision of VI. The optimization problem \eqref{eq:UCopt1} is therefore augmented as follows:
\begin{subequations}
\label{eq:UCopt2}
\begin{alignat}{3}
    & \min_\Phi && \sum_{t\in \mathcal{T}} \sum_{i\in \mathcal{I}} \big(C^{\mathrm{SU}}_{i}y_{it} + C_{i} p_{it}\big) +
    \sum_{t\in \mathcal{T}} \sum_{v\in \mathcal{V}} C^{\mathrm{VI}}_{v}h_{vt}  \label{eq:UC_obj_withInertia}\\
    & \;\mathrm{s.t.} \quad && \text{constraints} \;\; \eqref{eq:pow_bal_da}\text{-}\eqref{eq:vardef}, \nonumber\\
    & && k_{v t} = \frac{p_{v t}}{\sum_{v\in \mathcal{V}} \widebar{P}_{v}}, \; \forall v,t, \label{eq:VI_K}\\
    & && \ushort{P}_v \le p_{vt} \leq \widebar{P}_v, \; \forall v,t, \label{eq:pmax_bat}\\
	& && h_{i t} =  2H_i \widebar{P}_i u_{it}, \; \forall i,t, \label{eq:inertia_g} \\
	& && h_{v t} =  2H_v p_{v t}, \; \forall v,t \label{eq:inertia_v} \\
	& && M_t^\mathrm{SG} =  \sum_{i\in \mathcal{I}} 2 H_i k_{i t}, \; \forall t, \label{eq:inertia_g_tot} \\
	& && M_t^\mathrm{VI} =  \sum_{v\in \mathcal{V}} 2 H_v k_{v t}, \; \forall t,\label{eq:inertia_v_tot} \\
	& && M_{t} = \frac{M_t^\mathrm{SG} \sum_{i\in \mathcal{I}} \widebar{P}_i + M_t^\mathrm{VI} \sum_{v\in \mathcal{V}} \widebar{P}_{v}}{\sum_{i\in \mathcal{I}} \widebar{P}_i + \sum_{v\in \mathcal{V}} \widebar{P}_{v}}, \; \forall t, \label{eq:iner_M_withInertia}\\
	& && h_{it} \ge 0, \forall i,t; p_{vt},k_{vt},h_{vt} \ge 0, \forall v,t, \label{eq:vardef_new}
\end{alignat}
\end{subequations}
with $\Phi = \{p_{it},u_{it},y_{it},k_{it},h_{i t}, \forall i,t; \delta_{nt}, \forall n,t; p_{v t},k_{v t},h_{v t},\forall v,t;$ $M_{t},M_{t}^\mathrm{SG},M_{t}^\mathrm{VI},$ $\forall t\}$ being the set of optimization variables and $C^\mathrm{VI}_{v}$ reflecting the cost of VI for each wind generator $v\in\mathcal{V}$.

Expression \eqref{eq:VI_K} defines the power scaling factors $k_{v t}$ of wind farms providing VI\footnote{Note that the VI scaling factors $k_{v t}$ are conceptually different from the SG scaling factors $k_{i t}$ in \eqref{eq:iner_K}, since the time-variant aspect comes from the active power output $p_{v t}$ of the VI device and not the binary decision variable $u_{v t}$. This is justified by the fact that the VI devices are assumed to always be online (i.e., available for FFR provision) and hence $u_{v t}=1, \forall v,t$.} as a function of the respective battery power output $p_{vt}$, bounded by the upper and lower limits in \eqref{eq:pmax_bat}, whereas constraints \eqref{eq:inertia_g}-\eqref{eq:inertia_v} describe the individual inertia constants\footnote{Due to discrete nature of the objective function \eqref{eq:UC_obj_withInertia}, the inertia constants $h_{i t}$ and $h_{v t}$ are expressed in \SI{}{\mega\watt\second\squared}.} of each conventional and VI-providing generator respectively. The aggregate inertia contributions from synchronous and converter-interfaced generators are accounted for in \eqref{eq:inertia_g_tot}-\eqref{eq:inertia_v_tot}, \eqref{eq:iner_M_withInertia} defines the global system inertia as a weighted average of inertia gains provided by SGs and wind farms, while \eqref{eq:vardef_new} declares the new decision variables.

\vspace{-0.35cm}
\section{Economic Valuation of Inertia} \label{sec:EconVal}

\subsection{Value of Inertia for the TSO} \label{sec:InertiaCost}

The consideration of minimum inertia requirements for a TSO raises the question of how much these requirements will cost. To assess the economic value of inertia, the UC results can be used to determine the difference in system costs with and without inertia requirements, which would in turn provide an upper bound of the value of inertia for the TSO. 

As an illustrative example, Table~\ref{tab:UC_costs} shows the breakdown of UC costs for the modified IEEE RTS-96 system from \cite{modifiedRTS96}, for cases with and without inertia requirements. Several conclusions can be drawn from Table~\ref{tab:UC_costs}. First, it is clear that the inclusion of inertia requirements has a significant impact on the system cost, with the total commitment costs increasing by \SI{162100}{\sieuro}, which indicates that the newly committed SGs are needed solely for satisfying the inertia requirements. Finally, a large increase in generator dispatch costs reflects the energy produced by the SGs additionally committed for inertia purposes, thus suggesting that they are operating at their technical minimum.

\begin{table}[!t]
\renewcommand{\arraystretch}{1.2}
\caption{UC costs $\mathrm[\text{\euro}]$ with and without inertia requirements.}
\label{tab:UC_costs}
\noindent
\centering
    \begin{minipage}{\linewidth} 
    \renewcommand\footnoterule{\vspace*{-5pt}} 
    \begin{center}
        \begin{tabular}{ c || c | c | c }
            \toprule
            \textbf{Case} & \textbf{Total cost} & \textbf{Start-up cost} & \textbf{Energy cost} \\ 
            \cline{1-4}
            w/o inertia requirements & $407\,500$ & $600$ & $406\,900$\\
            w/ inertia requirements & $569\,600$ & $3\,500$ & $566\,100$\\
            \arrayrulecolor{black}\bottomrule
        \end{tabular}
        \end{center}
    \end{minipage}
    \vspace{-0.35cm}
\end{table}

It can be concluded from the aforementioned results that the difference in cost between these two cases of \SI{162100}{\sieuro} can be used as an indicator for defining the price at which the TSO values inertia. In other words, this sum represents the maximum amount the TSO would be willing to pay for procuring inertia from other units.

\vspace{-0.35cm}
\subsection{Inertia Substitution} \label{sec:InertiaSubstitution}

Authors in \cite{xu2016economic} and \cite{badesa2017economic} employ the method of adding ``free'' inertia to a frequency-constrained unit commitment and defining the decrease in total costs resulting from this addition as the value of inertia. Such concept displaces the more expensive inertia provided by synchronous generators with the cheap inertia and simultaneously preserves the total system-level inertia constant. While adding free inertia to the system is an effective tool for understanding the marginal inertia cost of the traditional SG portfolio, it does not reflect the cost of inertia provision for converter-interfaced generators. Building on a similar approach, the technique used in this work is to include the continuous inertia slack variable $M_{t}^+ \in \R_{\geq0}$ in the RoCoF constraint of the UC problem, with an associated cost $C^+$ present in the objective function. The RoCoF constraint \eqref{eq:rocof_cstr} is thus reformulated as $\frac{\dot{f}_\mathrm{lim}}{f_0} (M_{t}+M_{t}^+) \geq \lvert\Delta P_{t}\rvert,\,\forall t$, whereas the objective function \eqref{eq:UC_obj} is adapted accordingly:
\begin{equation}
    \begin{aligned}
    &\min_\Phi \,\, \sum_{t\in \mathcal{T}} \sum_{i\in \mathcal{I}} \big(C^{\mathrm{SU}}_{i}y_{it} + C_{i} p_{it}\big)\ +
    \sum_{t\in \mathcal{T}} C^+ M_{t}^+. \label{eq:UC_obj_mod}
    \end{aligned}
\end{equation}
It should be noted that the virtual inertia is not considered in this case, as it will later replace the slack variable. Using the described method, a certain cost $C^+$ for which the additional inertia is no longer profitable can be found, reflected by the optimizer deciding to commit conventional generators instead.


\begin{figure}[!b]
    \centering
    \vspace{-0.35cm}
    \scalebox{0.0485}{\includegraphics[]{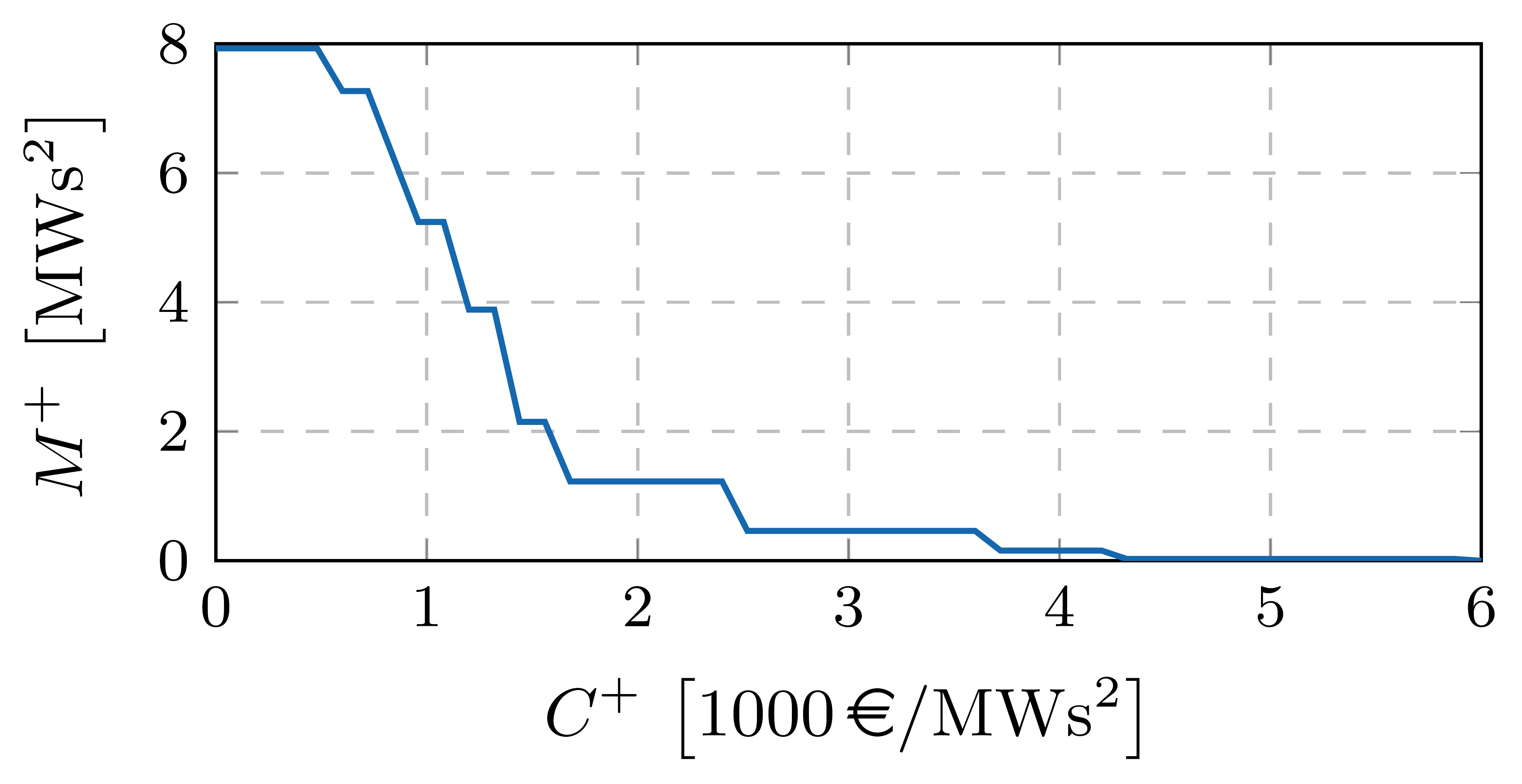}}
    \caption{Cost curve of additional inertia provision.}
    \label{fig:substitution_curve}
\end{figure}

For illustration purposes, the UC is ran on the previously described test case for a wide range of $C^+$ values. For each run, the amount of extra inertia purchased by the system is recorded, with the final results illustrated in Fig.~\ref{fig:substitution_curve}. Interestingly, the relationship between $M^+$ and $C^+$ is reminiscent of a demand curve and it can be interpreted in a similar fashion. The main difference compared to a traditional demand curve lies in the fact that the total amount of system inertia does not change for different values of $M^+$. More precisely, the RoCoF constraint ensures that the minimum amount of inertia needed for system stability is obtained, and including cheaper inertia does not increase the total inertia in the system but rather substitutes more expensive inertia provided by conventional synchronous generators instead. This effect is clearly reflected by distinctive steps in the extra inertia cost curve, with each step corresponding to an SG being displaced by a unit with cheaper inertia provision. Eventually, at $M^+\approx\SI{8}{\mega\watt\second\squared}$ the system does not purchase any more slack inertia despite it being free. This is justified by the fact that all expensive SGs, initially committed solely for providing inertia and meeting the frequency requirements, have now been decommissioned and replaced by addition of cheap inertia. Similarly, the extra inertia will not be procured at a cost higher than $6000\,\text{\euro}/\mathrm{MW\,s}^2$, since it is more expensive than having all necessary SGs come online.


\vspace{-0.35cm}
\section{Inertia Pricing \& Payment Methods} \label{sec:InertiaPricing}
In this section, three novel methods for pricing and payment of inertia are proposed. In particular, the prices are derived from the UC model by obtaining the dual variables of the respective constraints. Nonetheless, this procedure requires to relax the original Mixed-Integer Linear Program (MILP) form of the UC problem into a Linear Program (LP) whose dual variables are defined explicitly. For instance, the price of energy is determined as the dual of the power balance constraint, and there are three alternative price models for obtaining this dual \cite{gribik2007market}: (i) the restricted model; (ii) the dispatchable model; and (iii) the convex hull model.

In the restricted model, binary variables $u$ are set to the optimal value $u^*$ according to the solution of the preceding MILP and the optimization problem is solved again in the LP form. In contrast, the pricing in the dispatchable model is obtained by allowing binary variables to vary continuously in the range $0 \leq u \leq 1$. Finally, the convex hull model seeks to approximate the non-convex aggregate cost function with a convex hull function \cite{gribik2007market}. In this paper we will solely focus on the restricted model for price determination, given that the main goal of our work is to investigate alternative inertia pricing schemes and not to focus on the non-convex nature of the UC problem. Indeed, while the convex hull method sidesteps the non-convexity issues it is also very computationally intensive, whereas the dispatchable model introduces errors in the inertia calculation based on variable $u$.

\vspace{-0.35cm}
\subsection{Issues Pertaining to  Negative SG Profits} \label{sec:ProfitAnalysis}

The inclusion of inertia requirements in the UC forces a dispatch of additional SGs solely for the purpose of providing rotational inertia, which in turn preserves the frequency deviation within prescribed limits after a fault. Under such circumstances, the expensive SGs are started up and operated at their technical minimum to ensure grid stability. Nonetheless, with extra units coming online, the total generator profits are affected, which in turn calls for the establishment of proper reimbursement mechanisms to make whole those units which should be part of the efficient dispatch but are not supported by the current prices.


The increase in the number of online SGs can potentially lead to negative profits for units committed solely for inertia purposes, as they are dispatched out of the energy merit order. Indeed, the hourly price of energy is determined as the marginal cost of supplying the load. However, if more expensive generators come online to provide inertia by operating at the technical minimum, their marginal costs will not be reflected in the energy price. Hence, their fuel costs remain above the market price which eventually leads to negative revenue. In addition to the loss in the Energy-Only Market (EOM), they also face losses from the start-up costs. The issue of generators facing negative profits is important as there is no incentive for them to come online for inertia provision. The TSO must therefore create payment mechanisms to reimburse SGs and encourage cheaper and more efficient technologies to provide inertia.


\vspace{-0.35cm}
\subsection{Limitations of Obtaining Duals from RoCoF Constraint} \label{sec:RocofDual}

To obtain a price for inertia and create payments to units providing it, the constraint of interest is the upper bound on RoCoF imposed in \eqref{eq:rocof_cstr}. Due to non-convexity issues, obtaining a dual from such constraint is not straightforward. When considering SG units, the cost of providing inertia is linked to both the start-up costs and the fuel costs. The former costs depend on the binary commitment variables which complicates the extraction of the dual variables. Indeed, when only considering inertia provided by SGs, the RoCoF constraint is satisfied by committing a certain amount of generators. However, since the inertia comes in steps due to the binary nature of the problem, the constraint will be satisfied once the threshold is exceeded, but it will (almost) never be equal to the exact limit. Therefore, the dual of the RoCoF constraint would be zero in this particular case.

Therefore, the main issue to be addressed in order to extract a meaningful RoCoF dual variable is the discrete online/offline characteristic of SGs and the fact that they do not provide inertia in a continuous manner. Hence, a straightforward approach to obtain a RoCoF dual variable would be to set the SG inertia as a continuous variable and rewrite constraint \eqref{eq:iner_K} as
\begin{equation}
    0 \leq k_{i t} \leq \frac{\widebar{P}_{i} }{\sum_{i\in \mathcal{I}} \widebar{P}_{i}} u_{it},   \quad \forall i,t. \label{eq:iner_K_mod}
\end{equation}
The aggregate inertia remains equal to $M_{t} = \sum_{i\in \mathcal{I}} 2 H_i k_{i t}, \forall t$, and is therefore transformed into a continuous variable.

While such formulation enables SGs to artificially provide inertia in a continuous fashion, it still does not yield a meaningful RoCoF dual variable. More precisely, the problem lies in the fact that there is no direct cost associated with $M_t$ in the objective function \eqref{eq:UC_obj}. Duals represent how a change in the respective constraint would affect the value of the objective function. However, the cost of inertia provision from SGs is not explicitly included in the objective function, but it is in fact ``hidden'' in the start-up and fuel costs.

On the other hand, when considering virtual inertia in the UC, a non-zero dual variable can be obtained for RoCoF constraint. Nevertheless, it only takes the bid-in cost $C_v^\mathrm{VI}$ of converter-interfaced units into account, since the objective function \eqref{eq:UC_obj_withInertia} has an explicit cost term $C_v^\mathrm{VI}h_{vt}$ for purchasing VI. This however does not resolve the aforementioned issue pertaining to traditional SG units providing inertia, as their costs are still not captured in the existing price setting mechanism. Hence, alternative pricing methods must be developed instead.

\vspace{-0.35cm}
\subsection{Proposed Inertia Pricing and Payment Methods} \label{sec:PricingMethods}

\begin{figure}[!t]
    \centering
    \scalebox{0.675}{\includegraphics[]{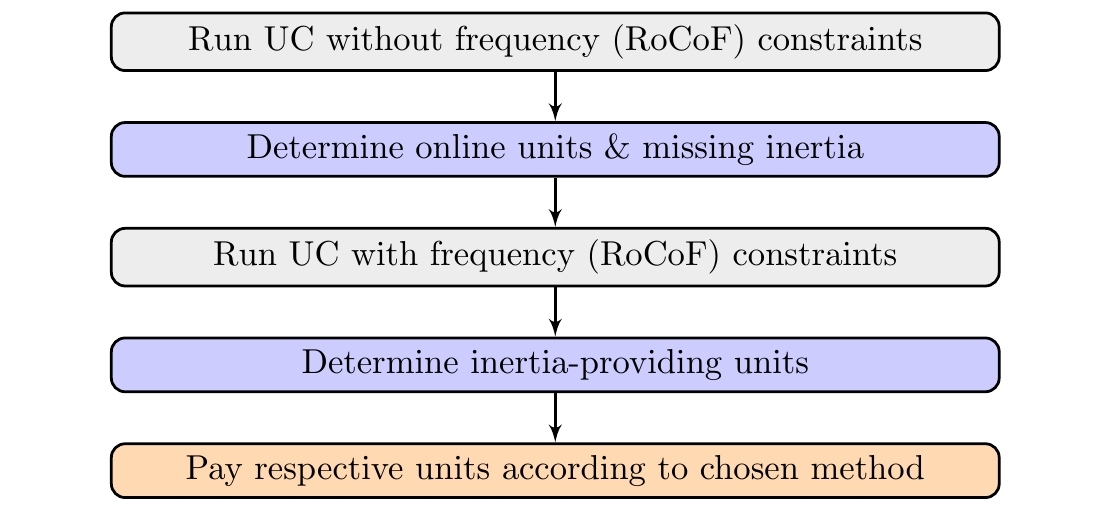}}
    \caption{Proposed methodology for determining inertia payments.}
    \label{fig:methodology}
    \vspace{-0.35cm}
\end{figure}

Three novel methods for pricing and payment of inertia for all participating units are developed in this study, all using the same two-step UC algorithm described in Fig.~\ref{fig:methodology}. The main purpose of this algorithm is to differentiate between the units being online for energy purposes and the ones committed additionally solely for inertia provision. This breakdown allows payment schemes to reimburse units in a fair manner according to their individual participation. Moreover, the proposed two-step process is analogous to the well-known two-step reliability commitment process \cite{gribik2007market}. The proposed pricing methods are described in detail in the remainder of this sections.


\subsubsection{Ex-Post Price Calculation}

The first pricing method is based on \textit{ex post} calculations of the inertia price. Motivated by the clearing price determination in an EOM, this method aims at setting the price of inertia provision as the cost of the most expensive unit supplying the inertia demand.

Let us assume obtaining a dual $\lambda^\mathrm{H}_t, \forall t$ from the RoCoF constraint \eqref{eq:rocof_cstr}, which would in this case represent the cost of procuring virtual inertia. In order to include the costs incurred to SGs providing rotational inertia into analysis, the following price structure is introduced: 
\begin{equation}
    \hat{\lambda}^\mathrm{H}_t = \max \big\{\lambda^\mathrm{H}_t, \lambda_t^\mathrm{SG} \big\}, \quad \forall t \in \mathcal{T}, \label{eq:m1_def}
\end{equation}
with $\hat{\lambda}^\mathrm{H}_t$ being the marginal value of inertia provision at hour $t$ and $\lambda_t^\mathrm{SG} = \max_{i\in\mathcal{M}^+} \left\{ \left(\max \big\{C_i-\mu_t, 0 \big\} \ p_{it} + C_i^\mathrm{SU}\right)u_{it}/h_{i t}\right\}$ determining the cost of the most expensive SG unit online (and hence providing inertia); $\mathcal{M}^+\in\N$ represents the set of additional SGs committed in the second UC run and $\mu_t$ is the dual of the power balance constraint corresponding to the energy price.

According to \eqref{eq:m1_def}, the price of inertia is defined as the cost of the most expensive unit (VI or SG) supplying inertia. The SG costs are obtained by combining the EOM losses, assuming the generator is in fact facing any, together with the start-up costs, and subsequently dividing it by the amount of inertia provided. Once the price is determined, units get reimbursed through Inertia Payments (IP) based on the previously determined price $\hat{\lambda}^\mathrm{H}_t$ and the quantity of inertia that they provide ($h_{it}$ for SGs and $h_{vt}$ for VI providers), scaled accordingly as follows:
\begin{subequations}
\label{eq:inertia_payments}
\begin{align}
    \mathrm{IP}_{it} &=  \hat{\lambda}^\mathrm{H}_t \frac{h_{it} \sum_{i\in \mathcal{I}} \widebar{P}_i}{\sum_{i\in \mathcal{I}} \widebar{P}_i + \sum_{v\in \mathcal{V}} \widebar{P}_{v}},   \quad \forall i,t, \\
    \mathrm{IP}_{vt} &=  \hat{\lambda}^\mathrm{H}_t \frac{h_{v t}\sum_{v\in \mathcal{V}} \widebar{P}_{v}}{\sum_{i\in \mathcal{I}} \widebar{P}_i + \sum_{v\in \mathcal{V}} \widebar{P}_{v}},   \quad \forall v,t.
\end{align}
\end{subequations}

\subsubsection{Utility Function}

In Section~\ref{sec:InertiaCost} it was shown that the financial resources the TSO would be willing to spend on inertia can be quantified. Indeed, for the given test system, this value corresponds to \SI{162100}{\sieuro} for procuring \SI{8}{\mega\watt\second\squared} and can be defined as the utility of inertia for the TSO, thus providing information on how to price the inertia in the system.

In the same vein, this method aims to obtain the relevant pricing information from the UC by subtracting the utility term $U^\mathrm{H} \sum_{t} M_t$ from the objective function \eqref{eq:UC_obj_withInertia}, where $U^\mathrm{H}$ denotes the utility of inertia for the entire scheduling horizon. The cumulative value of $U^\mathrm{H}$ is previously calculated by running the UC twice for certain inertia demand: once without frequency constraints (denoted by ``$\mathrm{w/o~FC}$'') and once with frequency constraints and no VI provision (denoted by ``$\mathrm{w/~FC}$''). Such procedure yields the lower and upper bound on system cost, respectively. The utility of inertia is then determined by
\begin{equation}
    U^\mathrm{H} = \frac{C_\mathrm{total}^\mathrm{w/~FC} -C_\mathrm{total}^\mathrm{w/o~FC}}{H_\mathrm{dem}},
    \label{eq:utility}
\end{equation}
with $C_\mathrm{total}$ referring to the total UC cost. It should be noted that by inertia demand $H_\mathrm{dem}$ we consider the amount of additional inertia needed to satisfy the frequency requirements after the first UC run. Since the inertia costs of both VI and SG units are explicitly included in the objective function, a price for inertia can be obtained as the dual variable $\lambda^\mathrm{H}_t$ of the RoCoF constraint, i.e., $\hat{\lambda}^\mathrm{H}_t\coloneqq\lambda^\mathrm{H}_t$. The units would subsequently be reimbursed for inertia provision through the same payment mechanism described in \eqref{eq:inertia_payments}.

\subsubsection{Uplift Payments}

The authors in \cite{EnergiforskReport} build on the concept of uplift payments or ``commitment ticket price'' for generators providing inertia, previously introduced in \cite{hogan2003minimum}. They define the additional payments to the units as
\begin{equation}
    \pi_{it} = C_i^\mathrm{SU} - \widebar{\nu}_{it} \widebar{P}_i + \ushort{\nu}_{it} \ushort{P}_i,  \quad \forall i,t,
    \label{eq:uplift}
\end{equation}
where $\widebar{\nu}_{it}, \forall i,t$ and $\ushort{\nu}_{it}, \forall i,t$ are the duals of the maximum and minimum unit generation constraints, respectively. The authors then suggest to include the product $\mathrm{IP}_{it} = \pi_{it}\,y_{it}$ to the amount paid to the generators. The main drawback of the suggested payment lies in the fact that SGs are only compensated during start-up hours, as the payment is only made when $y_{it}>0$. It implies that the generators facing losses due to energy prices being lower than their marginal costs are not reimbursed during operating hours. Furthermore, an additional pricing mechanism for VI units must be implemented. 

In order to resolve the above limitations, the proposition in this work is to alter the uplift payment such that $\pi_{it} = C_i^\mathrm{SU}y_{it} + \ushort{\nu}_{it} \ushort{P}_i$\footnote{The term $\widebar{\nu}_{it} \widebar{P}_i$ is omitted for simplicity, since the probability of all additionally committed SGs operating at their technical maximum is rather negligible. However, the aforementioned term can be also incorporated if needed. Note that the case studies and results presented in the remainder of the chapter are not affected by this term.}. In contrast to \eqref{eq:uplift}, the modified uplift method allows the generators to receive payments for start-up costs as well as for all time periods in which they face EOM losses, and not only during the start-up hour. Indeed, the start-up decision variable $y_{it}$ has been incorporated to reflect the incurred start-up costs, but SGs will receive additional compensation when operating at the technical minimum reflected by the expression $\ushort{\nu}_{it}P_i^\mathrm{min}$. Nevertheless, such uplift payment is not sufficient for accurately capturing the costs of VI units. Hence, the VI providers will be reimbursed according to the dual $\lambda^\mathrm{H}_t$ of the RoCoF constraint obtained from the optimization, which still reflects the cost of inertia coming from such units. This dual represents in fact the VI common clearing price used to reimburse all virtual inertia units. The proposed payment method thus distinguishes between a standard shadow price mechanism for VI units and a reimbursement method based on traditional uplift payments for conventional SG units.

\vspace{-0.15cm}
\section{Case Studies}
\label{sec:CaseStudy}


\subsection{Small Test System} \label{sec:SmallSystem}


In order to gain basic understanding of the proposed pricing methodologies as well as the economic interactions in a UC, we first consider a simple system comprising three conventional SGs, fixed amount of wind generation with constant power output, and three batteries offering VI capabilities. The most relevant unit parameters including their normalized inertia constants $H_i$ and $H_v$ are given in Table~\ref{tab:small_system_description}. Note that the cost of batteries is not included as their bid-in cost for inertia $\left(C_v^\mathrm{VI}\right)$ is the only cost considered and assumed to be determined by the operators of VI units. For cases where battery costs are not a varying parameter, we assume that all batteries bid-in at a fixed price of $C_v^\mathrm{VI}=50\,\text{\euro}/\mathrm{MW\,s}^2$. This allows for dispatching of VI without completely removing the need for SGs, thus making the comparison between payment methods more comprehensive.

We consider an 8-hour scheduling horizon with the load demand varying between $180\,\mathrm{MWh}$ and $200\,\mathrm{MWh}$, and a constant wind energy output of $150\,\mathrm{MWh}$ at each hour. The inertia demand is dictated by the requirement for withstanding a predefined change in active power of $\Delta P = (1,1,15,25,33,33,30,10)\!\times\!10^{-3}\,\mathrm{p.u.}$ throughout the scheduling horizon, with smaller values indirectly corresponding to a lower disturbance probability.

\begin{table}[!t]
\renewcommand{\arraystretch}{1.2}
\caption{Small test system description.}
\label{tab:small_system_description}
\noindent
\centering
    \begin{minipage}{\linewidth} 
    \renewcommand\footnoterule{\vspace*{-5pt}} 
    \begin{center}
    \resizebox{\columnwidth}{!}{%
        \begin{tabular}{ c || c | c | c | c}
            \toprule
            \textbf{Unit} & $\boldsymbol{\widebar{P}_{i/v}/\ushort{P}_{i/v}\,[\mathrm{MW}]}$ &  $\boldsymbol{C_i\,[\text{\textbf{\euro}}/\mathrm{MWh}]}$ & $\boldsymbol{C_i^\mathrm{SU}\,[\text{\textbf{\euro}}]}$ & $\boldsymbol{H_{i/v}\,[\mathrm{s}]}$ \\ 
            \cline{1-5}
            $G_1$ & $160/10$ & $10$ & $0$ & $4$ \\
            $G_2$ & $100/10$ & $12$ & $300$ & $4$ \\
            $G_3$ & $80/10$ & $11$ & $200$ & $4$ \\
            $B_1$ & $10/0$ & $-$ & $-$ & $10$ \\
            $B_2$ & $20/0$ & $-$ & $-$ & $10$ \\
            $B_3$ & $30/0$ & $-$ & $-$ & $10$ \\
            \arrayrulecolor{black}\bottomrule
        \end{tabular}
    }
        \end{center}
    \end{minipage}
    \vspace{-0.35cm}
\end{table}

\begin{table}[!b]
\renewcommand{\arraystretch}{1.2}
\vspace{-0.35cm}
\caption{Parameters of VI battery units in the large system.}
\label{tab:IEEE_RTS96_batteries}
\noindent
\centering
    \begin{minipage}{\linewidth} 
    \renewcommand\footnoterule{\vspace*{-5pt}} 
    \begin{center}
        \begin{tabular}{ c || c | c | c | c}
            \toprule
            \textbf{Unit} & \textbf{Origin} & $\boldsymbol{\widebar{P}_v\,[\mathrm{MW}]}$ & $\boldsymbol{H_v\,[\mathrm{s}]}$ & $\boldsymbol{C_v^\mathrm{VI}\,[\text{\textbf{\euro}}/\mathrm{MW\,s}^2]}$ \\
            \cline{1-5}
            $B_1$ & $B102$ & $70$ & $12$ & $10$ \\
            $B_2$ & $B114$ & $50$ & $12$ & $11$ \\
            $B_3$ & $B202$ & $100$ & $12$ & $12$ \\
            $B_4$ & $B212$ & $40$ & $12$ & $13$ \\
            \arrayrulecolor{black}\bottomrule
        \end{tabular}
        \end{center}
    \end{minipage}
\end{table}

\vspace{-0.35cm}
\subsection{Large Test System} \label{sec:LargeSystem}

For a more comprehensive test case, we focus on the IEEE RTS-96 system \cite{modifiedRTS96} previously introduced in Section~\ref{sec:InertiaCost}, modified into a low-inertia system for the purposes of this study by adding a large amount of wind generation. For the analysis of virtual inertia pricing, four VI-controlled batteries are included in the system, with relevant parameters provided in Table~\ref{tab:IEEE_RTS96_batteries}. The costs shown are the bid-in costs for virtual inertia, as it is assumed throughout this work that the battery operators make decisions on which markets to participate in and thus select the appropriate bid-in costs to ensure profitable operation. In our case the battery operators enter the inertia provision market with a specific bid-in cost and are dispatched only if it is financially viable for the system. The given costs are not necessarily kept constant throughout this study and will occasionally be changed for certain test cases.

The scheduling horizon is a typical 24-hour day-ahead time period. Demand and wind production profiles are taken from \cite{modifiedRTS96}. Moreover, the system is required to be able to maintain frequency within prescribed thresholds after a disturbance of $\Delta P = 0.1\,\mathrm{p.u.}$ at each hour of the day.

The following sections present the performance of each pricing method on both the small and the large test system, as well as with and without the inclusion of virtual inertia.

\vspace{-0.35cm}
\subsection{Method 1: Ex-Post Price Calculation}

\subsubsection{Small test system} 
Without enforcing constraints on inertia, only $G_1$ is needed to supply the demand. The total UC cost is \SI{3 360}{\sieuro} and $G_1$ makes zero profit since it is the only generator online and thus sets the EOM price (equal to $10\,\text{\euro}/\mathrm{MWh}$). With the addition of inertia requirements, units $G_2$ and $G_3$ come online for $3$ hours each, resulting in a \SI{590}{\sieuro} increase in UC cost. Due to EOM price remaining the same, both units face losses from the EOM and start-up costs, leading to negative profits of \SI{-360}{\sieuro} and \SI{-230}{\sieuro} for $G_2$ and $G_3$ respectively. Table~\ref{tab:M1_results_small_system} provides more insights regarding the revenue, costs, payments and profits for $G_2$ and $G_3$ over the 8-hour scheduling horizon.

\begin{table}[!t]
\renewcommand{\arraystretch}{1.2}
\caption{Breakdown of hourly UC costs [\euro] for \textit{Method~1}.}
\label{tab:M1_results_small_system}
\noindent
\centering
    \begin{minipage}{\linewidth} 
    \renewcommand\footnoterule{\vspace*{-5pt}} 
    \begin{center}
    \resizebox{\columnwidth}{!}{%
        \begin{tabular}{ c || l | c | c | c | c | c | c | c | c}
            \toprule
            \textbf{Unit} & \textbf{Cost} & $\boldsymbol{H_1}$ & $\boldsymbol{H_2}$ & $\boldsymbol{H_3}$ & $\boldsymbol{H_4}$ & $\boldsymbol{H_5}$ & $\boldsymbol{H_6}$ & $\boldsymbol{H_7}$ & $\boldsymbol{H_8}$ \\
            \cline{1-10}
            \multirow{4}{*}{$\boldsymbol{G_2}$} & EOM & $-$ & $-$ & $-$ & $-$ & $-20$ & $-20$ & $-20$ & $-$  \\
                                 & SU    & $-$ & $-$ & $-$ & $-$ & $-300$ & $0$ & $0$ & $-$     \\
                                 & IP    & $-$ & $-$ & $-$ & $-$ & $320$ & $20$ & $20$ & $-$    \\
                                 & Total & $-$ & $-$ & $-$ & $-$ & $0$ & $0$ & $0$ & $-$         \\
            \cline{1-10}
            \multirow{4}{*}{$\boldsymbol{G_3}$} & EOM & $-$ & $-$ & $-$ & $-10$ & $-10$ & $-10$ & $-$ & $-$   \\
                                 & SU    & $-$ & $-$ & $-$ & $-200$ & $0$ & $0$ & $-$ & $-$      \\
                                 & IP    & $-$ & $-$ & $-$ & $210$ & $65.6$ & $39.5$ & $-$ & $-$  \\
                                 & Total & $-$ & $-$ & $-$ & $0$ & $55.6$ & $29.5$ & $-$ & $-$     \\
            \arrayrulecolor{black}\bottomrule
        \end{tabular}
        }
        \end{center}
    \end{minipage}
    \vspace{-0.35cm}
\end{table}

\begin{table}[!b]
\renewcommand{\arraystretch}{0.8}
\vspace{-0.35cm}
\caption{Hourly SG profits [\euro] under \textit{Method~1} pricing scheme.}
\label{tab:SG_profits_M1}
\noindent
\centering
\begin{minipage}{1\linewidth}
\renewcommand\footnoterule{\vspace*{-5pt}}
\begin{center}
    \huge
\resizebox{\columnwidth}{!}{%
\begin{tabular}{ c"c"c"c"c"c"c"c"c"c"c"c"c"c"c"c"c"c"c"c"c}
    \thickhline
    \cellcolor{black!10}$\Huge\boldsymbol{H\text{\textbackslash}G}$& $\Huge\boldsymbol{G_1}$ & $\Huge\boldsymbol{G_2}$ & $\Huge\boldsymbol{G_3}$ & $\Huge\boldsymbol{G_4}$ & $\Huge\boldsymbol{G_5}$ & $\Huge\boldsymbol{G_6}$ & $\Huge\boldsymbol{G_7}$ & $\Huge\boldsymbol{G_8}$ & $\Huge\boldsymbol{G_9}$ & $\Huge\boldsymbol{G_{10}}$ & $\Huge\boldsymbol{G_{11}}$ & $\Huge\boldsymbol{G_{12}}$ & $\Huge\boldsymbol{G_{13}}$ & $\Huge\boldsymbol{G_{14}}$ & $\Huge\boldsymbol{G_{15}}$ & $\Huge\boldsymbol{G_{16}}$ & $\Huge\boldsymbol{G_{17}}$ & $\Huge\boldsymbol{G_{18}}$ & $\Huge\boldsymbol{G_{19}}$ & $\Huge\boldsymbol{G_{20}}$ \\
    \thickhline
$\Huge\boldsymbol{H_1}$ & 	\cellcolor{myGreen!25}915 & 	\cellcolor{myGreen!25}915 & 	\cellcolor{myGreen!25}926 & 	\cellcolor{myGreen!25}1164 & 	\cellcolor{myGreen!25}891 & 	\cellcolor{white}0 & 	\cellcolor{white}0 & 	\cellcolor{white}0 & 	\cellcolor{myGreen!25}1531 & 	\cellcolor{white}0 & 	\cellcolor{myGreen!25}915 & 	\cellcolor{myGreen!25}915 & 	\cellcolor{myGreen!25}926 & 	\cellcolor{myGreen!25}1164 & 	\cellcolor{myGreen!25}891 & 	\cellcolor{white}0 & 	\cellcolor{white}0 & 	\cellcolor{white}0 & 	\cellcolor{myGreen!25}1531 & 	\cellcolor{myGreen!25}959	\\ \thickhline
$\Huge\boldsymbol{H_2}$ & 	\cellcolor{myGreen!25}211 & 	\cellcolor{myGreen!25}211 & 	\cellcolor{myGreen!25}214 & 	\cellcolor{white}0 & 	\cellcolor{myGreen!25}346 & 	\cellcolor{white}0 & 	\cellcolor{white}0 & 	\cellcolor{white}0 & 	\cellcolor{myGreen!25}394 & 	\cellcolor{myGreen!25}73 & 	\cellcolor{myGreen!25}211 & 	\cellcolor{myGreen!25}211 & 	\cellcolor{myGreen!25}214 & 	\cellcolor{white}0 & 	\cellcolor{myGreen!25}346 & 	\cellcolor{white}0 & 	\cellcolor{white}0 & 	\cellcolor{white}0 & 	\cellcolor{myGreen!25}394 & 	\cellcolor{myGreen!25}504	\\ \thickhline
$\Huge\boldsymbol{H_3}$ & 	\cellcolor{myGreen!25}211 & 	\cellcolor{myGreen!25}211 & 	\cellcolor{myGreen!25}214 & 	\cellcolor{white}0 & 	\cellcolor{myGreen!25}346 & 	\cellcolor{white}0 & 	\cellcolor{white}0 & 	\cellcolor{white}0 & 	\cellcolor{myGreen!25}394 & 	\cellcolor{myGreen!25}73 & 	\cellcolor{myGreen!25}211 & 	\cellcolor{myGreen!25}211 & 	\cellcolor{myGreen!25}214 & 	\cellcolor{white}0 & 	\cellcolor{myGreen!25}346 & 	\cellcolor{white}0 & 	\cellcolor{white}0 & 	\cellcolor{white}0 & 	\cellcolor{myGreen!25}394 & 	\cellcolor{myGreen!25}504	\\ \thickhline
$\Huge\boldsymbol{H_4}$ & 	\cellcolor{myGreen!25}211 & 	\cellcolor{myGreen!25}211 & 	\cellcolor{myGreen!25}214 & 	\cellcolor{white}0 & 	\cellcolor{myGreen!25}346 & 	\cellcolor{white}0 & 	\cellcolor{white}0 & 	\cellcolor{white}0 & 	\cellcolor{myGreen!25}394 & 	\cellcolor{myGreen!25}73 & 	\cellcolor{myGreen!25}211 & 	\cellcolor{myGreen!25}211 & 	\cellcolor{myGreen!25}214 & 	\cellcolor{white}0 & 	\cellcolor{myGreen!25}346 & 	\cellcolor{white}0 & 	\cellcolor{white}0 & 	\cellcolor{white}0 & 	\cellcolor{myGreen!25}394 & 	\cellcolor{myGreen!25}504	\\ \thickhline
$\Huge\boldsymbol{H_5}$ & 	\cellcolor{myGreen!25}211 & 	\cellcolor{myGreen!25}211 & 	\cellcolor{myGreen!25}214 & 	\cellcolor{white}0 & 	\cellcolor{myGreen!25}346 & 	\cellcolor{white}0 & 	\cellcolor{white}0 & 	\cellcolor{white}0 & 	\cellcolor{myGreen!25}394 & 	\cellcolor{myGreen!25}73 & 	\cellcolor{myGreen!25}211 & 	\cellcolor{myGreen!25}211 & 	\cellcolor{myGreen!25}214 & 	\cellcolor{white}0 & 	\cellcolor{myGreen!25}346 & 	\cellcolor{white}0 & 	\cellcolor{white}0 & 	\cellcolor{white}0 & 	\cellcolor{myGreen!25}394 & 	\cellcolor{myGreen!25}504	\\ \thickhline
$\Huge\boldsymbol{H_6}$ & 	\cellcolor{myGreen!25}211 & 	\cellcolor{myGreen!25}211 & 	\cellcolor{myGreen!25}214 & 	\cellcolor{white}0 & 	\cellcolor{myGreen!25}346 & 	\cellcolor{white}0 & 	\cellcolor{white}0 & 	\cellcolor{white}0 & 	\cellcolor{myGreen!25}394 & 	\cellcolor{myGreen!25}73 & 	\cellcolor{myGreen!25}211 & 	\cellcolor{myGreen!25}211 & 	\cellcolor{myGreen!25}214 & 	\cellcolor{white}0 & 	\cellcolor{myGreen!25}346 & 	\cellcolor{white}0 & 	\cellcolor{white}0 & 	\cellcolor{white}0 & 	\cellcolor{myGreen!25}394 & 	\cellcolor{myGreen!25}504	\\ \thickhline
$\Huge\boldsymbol{H_7}$ & 	\cellcolor{myGreen!25}75 & 	\cellcolor{myGreen!25}75 & 	\cellcolor{myGreen!25}214 & 	\cellcolor{red!25}-161 & 	\cellcolor{myGreen!25}116 & 	\cellcolor{white}0 & 	\cellcolor{myGreen!25}3966 & 	\cellcolor{myGreen!25}3966 & 	\cellcolor{myGreen!25}69 & 	\cellcolor{white}0 & 	\cellcolor{myGreen!25}75 & 	\cellcolor{myGreen!25}75 & 	\cellcolor{red!25}-46 & 	\cellcolor{red!25}-96 & 	\cellcolor{myGreen!25}116 & 	\cellcolor{white}0 & 	\cellcolor{myGreen!25}3966 & 	\cellcolor{myGreen!25}3966 & 	\cellcolor{myGreen!25}69 & 	\cellcolor{myGreen!25}171	\\ \thickhline
$\Huge\boldsymbol{H_8}$ & 	\cellcolor{red!25}-143 & 	\cellcolor{red!25}-143 & 	\cellcolor{white}0 & 	\cellcolor{myGreen!25}338 & 	\cellcolor{myGreen!25}94 & 	\cellcolor{white}0 & 	\cellcolor{myGreen!25}4424 & 	\cellcolor{myGreen!25}4424 & 	\cellcolor{myGreen!25}39 & 	\cellcolor{white}0 & 	\cellcolor{myGreen!25}65 & 	\cellcolor{myGreen!25}65 & 	\cellcolor{white}0 & 	\cellcolor{myGreen!25}274 & 	\cellcolor{myGreen!25}94 & 	\cellcolor{white}0 & 	\cellcolor{myGreen!25}4424 & 	\cellcolor{myGreen!25}4424 & 	\cellcolor{myGreen!25}39 & 	\cellcolor{myGreen!25}145	\\ \thickhline
$\Huge\boldsymbol{H_9}$ & 	\cellcolor{red!25}-143 & 	\cellcolor{red!25}-143 & 	\cellcolor{white}0 & 	\cellcolor{myGreen!25}338 & 	\cellcolor{myGreen!25}94 & 	\cellcolor{white}0 & 	\cellcolor{myGreen!25}4424 & 	\cellcolor{myGreen!25}4424 & 	\cellcolor{myGreen!25}39 & 	\cellcolor{white}0 & 	\cellcolor{myGreen!25}65 & 	\cellcolor{myGreen!25}65 & 	\cellcolor{white}0 & 	\cellcolor{myGreen!25}338 & 	\cellcolor{myGreen!25}94 & 	\cellcolor{white}0 & 	\cellcolor{myGreen!25}4424 & 	\cellcolor{myGreen!25}4424 & 	\cellcolor{myGreen!25}39 & 	\cellcolor{myGreen!25}145	\\ \thickhline
$\Huge\boldsymbol{H_{10}}$ & 	\cellcolor{red!25}-143 & 	\cellcolor{red!25}-143 & 	\cellcolor{white}0 & 	\cellcolor{myGreen!25}338 & 	\cellcolor{myGreen!25}94 & 	\cellcolor{white}0 & 	\cellcolor{myGreen!25}4424 & 	\cellcolor{myGreen!25}4424 & 	\cellcolor{myGreen!25}39 & 	\cellcolor{white}0 & 	\cellcolor{myGreen!25}65 & 	\cellcolor{myGreen!25}65 & 	\cellcolor{white}0 & 	\cellcolor{myGreen!25}338 & 	\cellcolor{myGreen!25}94 & 	\cellcolor{white}0 & 	\cellcolor{myGreen!25}4424 & 	\cellcolor{myGreen!25}4424 & 	\cellcolor{myGreen!25}39 & 	\cellcolor{myGreen!25}145	\\ \thickhline
$\Huge\boldsymbol{H_{11}}$ & 	\cellcolor{myGreen!25}65 & 	\cellcolor{myGreen!25}65 & 	\cellcolor{white}0 & 	\cellcolor{myGreen!25}338 & 	\cellcolor{myGreen!25}94 & 	\cellcolor{white}0 & 	\cellcolor{myGreen!25}4424 & 	\cellcolor{myGreen!25}4424 & 	\cellcolor{myGreen!25}39 & 	\cellcolor{white}0 & 	\cellcolor{myGreen!25}65 & 	\cellcolor{myGreen!25}65 & 	\cellcolor{white}0 & 	\cellcolor{myGreen!25}338 & 	\cellcolor{myGreen!25}94 & 	\cellcolor{white}0 & 	\cellcolor{myGreen!25}4424 & 	\cellcolor{myGreen!25}4424 & 	\cellcolor{myGreen!25}39 & 	\cellcolor{myGreen!25}145	\\ \thickhline
$\Huge\boldsymbol{H_{12}}$ & 	\cellcolor{myGreen!25}70 & 	\cellcolor{myGreen!25}70 & 	\cellcolor{red!25}-23 & 	\cellcolor{white}0 & 	\cellcolor{myGreen!25}105 & 	\cellcolor{white}0 & 	\cellcolor{myGreen!25}4195 & 	\cellcolor{myGreen!25}4195 & 	\cellcolor{myGreen!25}54 & 	\cellcolor{white}0 & 	\cellcolor{myGreen!25}70 & 	\cellcolor{myGreen!25}70 & 	\cellcolor{myGreen!25}217 & 	\cellcolor{white}0 & 	\cellcolor{myGreen!25}105 & 	\cellcolor{white}0 & 	\cellcolor{myGreen!25}4195 & 	\cellcolor{myGreen!25}4195 & 	\cellcolor{myGreen!25}54 & 	\cellcolor{myGreen!25}158	\\ \thickhline
$\Huge\boldsymbol{H_{13}}$ & 	\cellcolor{myGreen!25}70 & 	\cellcolor{myGreen!25}70 & 	\cellcolor{red!25}-23 & 	\cellcolor{white}0 & 	\cellcolor{myGreen!25}105 & 	\cellcolor{white}0 & 	\cellcolor{myGreen!25}4195 & 	\cellcolor{myGreen!25}4195 & 	\cellcolor{myGreen!25}54 & 	\cellcolor{white}0 & 	\cellcolor{myGreen!25}70 & 	\cellcolor{myGreen!25}70 & 	\cellcolor{myGreen!25}217 & 	\cellcolor{white}0 & 	\cellcolor{myGreen!25}105 & 	\cellcolor{white}0 & 	\cellcolor{myGreen!25}4195 & 	\cellcolor{myGreen!25}4195 & 	\cellcolor{myGreen!25}54 & 	\cellcolor{myGreen!25}158	\\ \thickhline
$\Huge\boldsymbol{H_{14}}$ & 	\cellcolor{myGreen!25}166 & 	\cellcolor{myGreen!25}166 & 	\cellcolor{red!25}-442 & 	\cellcolor{red!25}-1091 & 	\cellcolor{myGreen!25}313 & 	\cellcolor{white}0 & 	\cellcolor{white}0 & 	\cellcolor{white}0 & 	\cellcolor{myGreen!25}325 & 	\cellcolor{white}0 & 	\cellcolor{myGreen!25}166 & 	\cellcolor{myGreen!25}166 & 	\cellcolor{myGreen!25}167 & 	\cellcolor{red!25}-1091 & 	\cellcolor{myGreen!25}313 & 	\cellcolor{white}0 & 	\cellcolor{white}0 & 	\cellcolor{white}0 & 	\cellcolor{myGreen!25}325 & 	\cellcolor{myGreen!25}401	\\ \thickhline
$\Huge\boldsymbol{H_{15}}$ & 	\cellcolor{myGreen!25}211 & 	\cellcolor{myGreen!25}211 & 	\cellcolor{red!25}-442 & 	\cellcolor{white}0 & 	\cellcolor{myGreen!25}346 & 	\cellcolor{white}0 & 	\cellcolor{white}0 & 	\cellcolor{white}0 & 	\cellcolor{myGreen!25}394 & 	\cellcolor{myGreen!25}73 & 	\cellcolor{myGreen!25}211 & 	\cellcolor{myGreen!25}211 & 	\cellcolor{myGreen!25}214 & 	\cellcolor{white}0 & 	\cellcolor{myGreen!25}346 & 	\cellcolor{white}0 & 	\cellcolor{white}0 & 	\cellcolor{white}0 & 	\cellcolor{myGreen!25}394 & 	\cellcolor{myGreen!25}504	\\ \thickhline
$\Huge\boldsymbol{H_{16}}$ & 	\cellcolor{myGreen!25}211 & 	\cellcolor{myGreen!25}211 & 	\cellcolor{red!25}-442 & 	\cellcolor{white}0 & 	\cellcolor{myGreen!25}346 & 	\cellcolor{white}0 & 	\cellcolor{white}0 & 	\cellcolor{white}0 & 	\cellcolor{myGreen!25}394 & 	\cellcolor{myGreen!25}73 & 	\cellcolor{myGreen!25}211 & 	\cellcolor{myGreen!25}211 & 	\cellcolor{myGreen!25}214 & 	\cellcolor{white}0 & 	\cellcolor{myGreen!25}346 & 	\cellcolor{white}0 & 	\cellcolor{white}0 & 	\cellcolor{white}0 & 	\cellcolor{myGreen!25}394 & 	\cellcolor{myGreen!25}504	\\ \thickhline
$\Huge\boldsymbol{H_{17}}$ & 	\cellcolor{myGreen!25}211 & 	\cellcolor{myGreen!25}211 & 	\cellcolor{red!25}-442 & 	\cellcolor{white}0 & 	\cellcolor{myGreen!25}346 & 	\cellcolor{white}0 & 	\cellcolor{white}0 & 	\cellcolor{white}0 & 	\cellcolor{myGreen!25}394 & 	\cellcolor{myGreen!25}73 & 	\cellcolor{myGreen!25}211 & 	\cellcolor{myGreen!25}211 & 	\cellcolor{myGreen!25}214 & 	\cellcolor{white}0 & 	\cellcolor{myGreen!25}346 & 	\cellcolor{white}0 & 	\cellcolor{white}0 & 	\cellcolor{white}0 & 	\cellcolor{myGreen!25}394 & 	\cellcolor{myGreen!25}504	\\ \thickhline
$\Huge\boldsymbol{H_{18}}$ & 	\cellcolor{myGreen!25}211 & 	\cellcolor{myGreen!25}211 & 	\cellcolor{myGreen!25}214 & 	\cellcolor{white}0 & 	\cellcolor{myGreen!25}346 & 	\cellcolor{white}0 & 	\cellcolor{white}0 & 	\cellcolor{white}0 & 	\cellcolor{myGreen!25}394 & 	\cellcolor{myGreen!25}73 & 	\cellcolor{myGreen!25}211 & 	\cellcolor{myGreen!25}211 & 	\cellcolor{myGreen!25}214 & 	\cellcolor{white}0 & 	\cellcolor{myGreen!25}346 & 	\cellcolor{white}0 & 	\cellcolor{white}0 & 	\cellcolor{white}0 & 	\cellcolor{myGreen!25}394 & 	\cellcolor{myGreen!25}504	\\ \thickhline
$\Huge\boldsymbol{H_{19}}$ & 	\cellcolor{myGreen!25}211 & 	\cellcolor{myGreen!25}211 & 	\cellcolor{myGreen!25}214 & 	\cellcolor{white}0 & 	\cellcolor{myGreen!25}346 & 	\cellcolor{white}0 & 	\cellcolor{white}0 & 	\cellcolor{white}0 & 	\cellcolor{myGreen!25}394 & 	\cellcolor{myGreen!25}73 & 	\cellcolor{myGreen!25}211 & 	\cellcolor{myGreen!25}211 & 	\cellcolor{myGreen!25}214 & 	\cellcolor{white}0 & 	\cellcolor{myGreen!25}346 & 	\cellcolor{white}0 & 	\cellcolor{white}0 & 	\cellcolor{white}0 & 	\cellcolor{myGreen!25}394 & 	\cellcolor{myGreen!25}504	\\ \thickhline
$\Huge\boldsymbol{H_{20}}$ & 	\cellcolor{myGreen!25}464 & 	\cellcolor{myGreen!25}464 & 	\cellcolor{myGreen!25}377 & 	\cellcolor{white}0 & 	\cellcolor{myGreen!25}631 & 	\cellcolor{white}0 & 	\cellcolor{red!25}-723 & 	\cellcolor{red!25}-723 & 	\cellcolor{myGreen!25}872 & 	\cellcolor{myGreen!25}378 & 	\cellcolor{myGreen!25}464 & 	\cellcolor{myGreen!25}464 & 	\cellcolor{myGreen!25}377 & 	\cellcolor{white}0 & 	\cellcolor{myGreen!25}631 & 	\cellcolor{white}0 & 	\cellcolor{red!25}-723 & 	\cellcolor{red!25}-723 & 	\cellcolor{myGreen!25}872 & 	\cellcolor{myGreen!25}1107	\\ \thickhline
$\Huge\boldsymbol{H_{21}}$ & 	\cellcolor{myGreen!25}464 & 	\cellcolor{myGreen!25}464 & 	\cellcolor{myGreen!25}377 & 	\cellcolor{white}0 & 	\cellcolor{myGreen!25}631 & 	\cellcolor{white}0 & 	\cellcolor{red!25}-723 & 	\cellcolor{red!25}-723 & 	\cellcolor{myGreen!25}872 & 	\cellcolor{myGreen!25}378 & 	\cellcolor{myGreen!25}464 & 	\cellcolor{myGreen!25}464 & 	\cellcolor{myGreen!25}377 & 	\cellcolor{white}0 & 	\cellcolor{myGreen!25}631 & 	\cellcolor{white}0 & 	\cellcolor{red!25}-723 & 	\cellcolor{red!25}-723 & 	\cellcolor{myGreen!25}872 & 	\cellcolor{myGreen!25}1107	\\ \thickhline
$\Huge\boldsymbol{H_{22}}$ & 	\cellcolor{myGreen!25}464 & 	\cellcolor{myGreen!25}464 & 	\cellcolor{myGreen!25}377 & 	\cellcolor{white}0 & 	\cellcolor{myGreen!25}631 & 	\cellcolor{white}0 & 	\cellcolor{red!25}-723 & 	\cellcolor{red!25}-723 & 	\cellcolor{myGreen!25}872 & 	\cellcolor{myGreen!25}378 & 	\cellcolor{myGreen!25}464 & 	\cellcolor{myGreen!25}464 & 	\cellcolor{myGreen!25}377 & 	\cellcolor{white}0 & 	\cellcolor{myGreen!25}631 & 	\cellcolor{white}0 & 	\cellcolor{red!25}-723 & 	\cellcolor{red!25}-723 & 	\cellcolor{myGreen!25}872 & 	\cellcolor{myGreen!25}1107	\\ \thickhline
$\Huge\boldsymbol{H_{23}}$ & 	\cellcolor{myGreen!25}464 & 	\cellcolor{myGreen!25}464 & 	\cellcolor{myGreen!25}377 & 	\cellcolor{white}0 & 	\cellcolor{myGreen!25}631 & 	\cellcolor{white}0 & 	\cellcolor{red!25}-723 & 	\cellcolor{red!25}-723 & 	\cellcolor{myGreen!25}872 & 	\cellcolor{myGreen!25}378 & 	\cellcolor{myGreen!25}464 & 	\cellcolor{myGreen!25}464 & 	\cellcolor{myGreen!25}377 & 	\cellcolor{white}0 & 	\cellcolor{myGreen!25}631 & 	\cellcolor{white}0 & 	\cellcolor{red!25}-723 & 	\cellcolor{red!25}-723 & 	\cellcolor{myGreen!25}872 & 	\cellcolor{myGreen!25}1107	\\ \thickhline
$\Huge\boldsymbol{H_{24}}$ & 	\cellcolor{myGreen!25}464 & 	\cellcolor{myGreen!25}464 & 	\cellcolor{myGreen!25}377 & 	\cellcolor{white}0 & 	\cellcolor{myGreen!25}631 & 	\cellcolor{white}0 & 	\cellcolor{red!25}-723 & 	\cellcolor{red!25}-723 & 	\cellcolor{myGreen!25}872 & 	\cellcolor{myGreen!25}378 & 	\cellcolor{myGreen!25}464 & 	\cellcolor{myGreen!25}464 & 	\cellcolor{myGreen!25}377 & 	\cellcolor{white}0 & 	\cellcolor{myGreen!25}631 & 	\cellcolor{white}0 & 	\cellcolor{red!25}-723 & 	\cellcolor{red!25}-723 & 	\cellcolor{myGreen!25}872 & 	\cellcolor{myGreen!25}1107	\\
    \arrayrulecolor{black}\thickhline
\end{tabular}
 }
\end{center}
\end{minipage}
\end{table}

At hour $4$, $G_3$ turns on to provide inertia, and its start-up and EOM losses are covered by the inertia payments. Moreover, it turns on before $G_2$ due to its lower operating and start-up costs. Nevertheless, in the next hour $G_2$ also comes online, thus becoming the most expensive unit providing inertia. Inertia payments therefore cover its losses, resulting in $G_3$ making positive profit as it is no longer the most expensive unit providing inertia. Similar situation occurs at hour 6. Finally, in the next hour $G_3$ turns off and inertia payments cover the losses of generator $G_2$. Overall, $G_2$ and $G_3$ now have respective profits of \SI{0}{\sieuro} and \SI{85.1}{\sieuro}. The proposed payment method is thus successful in resolving negative profit issues and ensures that all units receive a common clearing inertia price.

\subsubsection{Large test system} 

\begin{figure}[!t]
    \centering
    \scalebox{0.0565}{\includegraphics[]{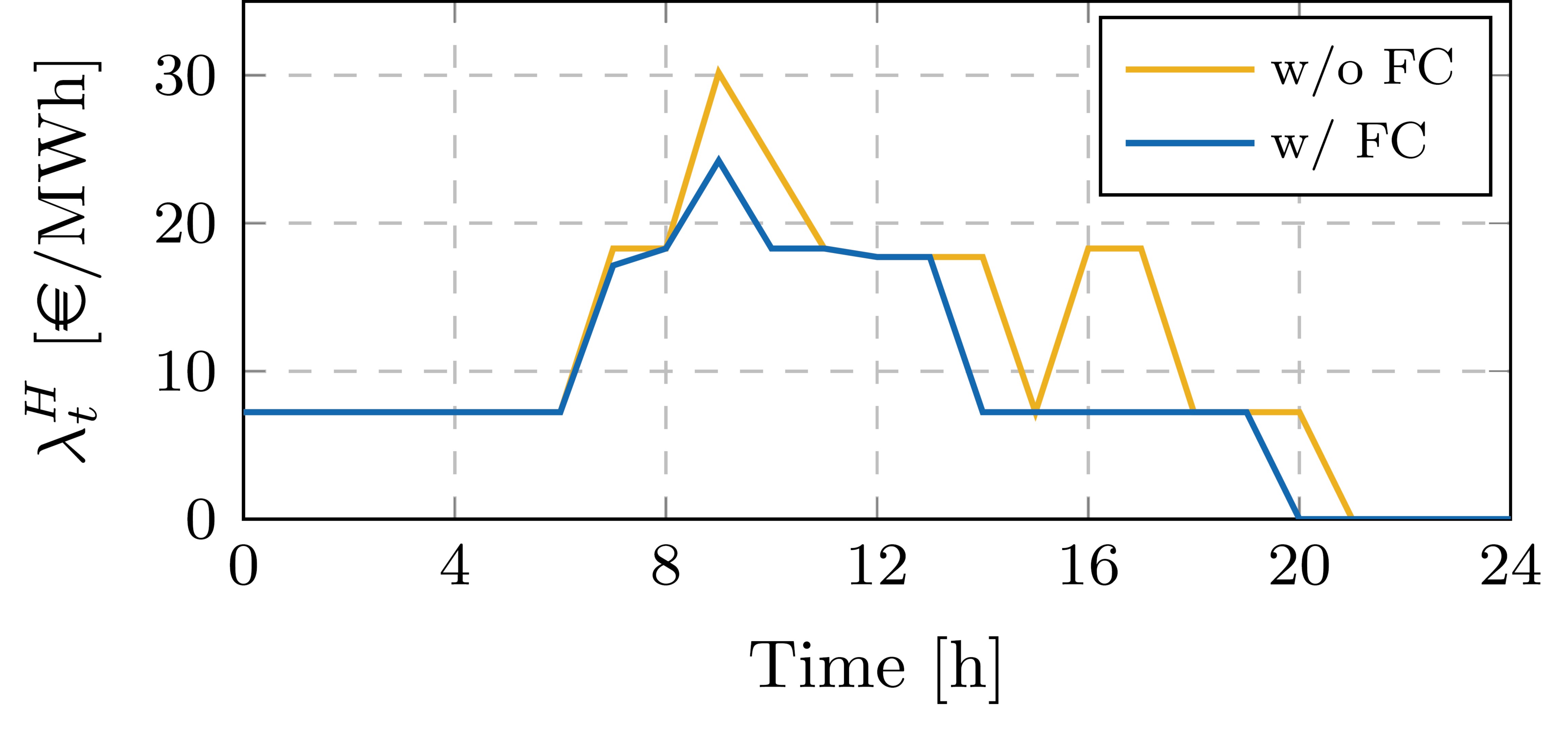}}
    \caption{EOM price with and without frequency constraints.}
    \label{fig:eom_comp}
    \vspace{-0.35cm}
\end{figure}

We investigate the performance of \textit{Method~1} on a larger system by analyzing the hourly profits of all $20$ generators given in Table~\ref{tab:SG_profits_M1}, with green and red fields indicating positive and negative hourly profits, respectively. The effectiveness of the payment scheme is clearly reflected in the vast majority of SG units achieving positive profits for most of the hours. The generators facing negative profits at some instances are the units participating in the EOM that suffer losses due to a decrease in energy price. More precisely, the inertia requirements and changed commitment of the system result in lower EOM prices and therefore reduced profits for units (see Fig.~\ref{fig:eom_comp}). This decrease in energy price is caused by the generators having to reduce their energy output to accommodate the extra energy in the system due to additional generators, mostly operating at their technical minimum, coming online. At certain times, this decrease of energy output is large enough to reduce the  price.


\begin{figure}[!b] 
	\centering
	\vspace{-0.35cm}
	\scalebox{0.0485}{\includegraphics[]{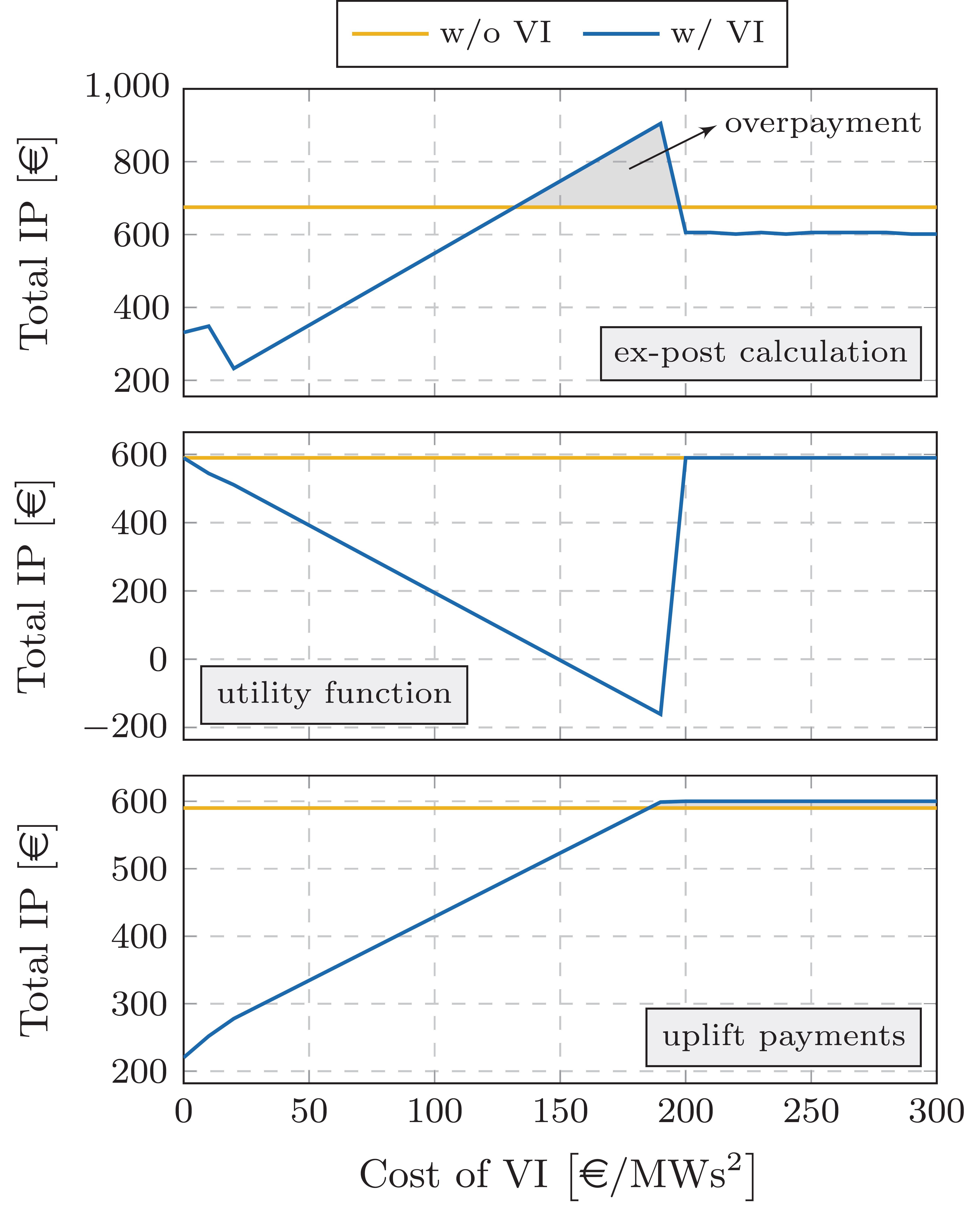}}
	\caption{Impact of virtual inertia on total inertia payments under different payment mechanisms: (i) ex-post calculation; (ii) utility function; (iii) uplift payments.}
	\label{fig:VI_payments}
\end{figure}

\subsubsection{Inclusion of Virtual Inertia}
The proposed method performs well with the addition of virtual inertia, but issues appear once the cost of VI provision becomes too expensive for the system. As an illustrative example, let us consider the small test system and assume a simple case where all VI units bid-in at the same cost for inertia. According to \eqref{eq:m1_def}, as this cost increases so do the corresponding inertia payments. However, a price threshold exists at which these payments exceed the corresponding amount for a respective case without any VI units in the system, as indicated in Fig.~\ref{fig:VI_payments}. While counter-intuitive at first, such effect can be justified by the fact that at certain hours the virtual inertia is cheaper than the start-up costs of a generator, and the TSO would rather employ VI instead of committing a conventional SG. However, the VI costs lead to higher duals in the subsequent hours and therefore higher total payments for inertia.

\vspace{-0.35cm}
\subsection{Method 2: Utility Function }
\subsubsection{Small test system} Similar to the analysis of \textit{Method~1}, Table~\ref{tab:M2_results_small_system} showcases the revenue, costs, payments and profits for generators $2$ and $3$ in the small test system under the payment scheme of \textit{Method~2}. For the entire scheduling horizon the total profits of $G_2$ and $G_3$ are \SI{-3}{\sieuro} and \SI{3}{\sieuro} respectively. Despite having opposite profits, the issue of larger negative profits due to inertia requirements is resolved by the proposed payment method. Indeed, if the opposite profits were to be redistributed \textit{ex-post}, then both $G_2$ and $G_3$ would break even and thus be indifferent to being committed for inertia provision. This symmetry in profits for SGs coming online for inertia provision is a specific characteristics of the utility function method and will be elaborated in more detail in the remainder of this section.

\begin{table}[!t]
\renewcommand{\arraystretch}{1.2}
\caption{Breakdown of hourly UC costs [\euro] for \textit{Method~2}.}
\label{tab:M2_results_small_system}
\noindent
\centering
    \begin{minipage}{\linewidth} 
    \renewcommand\footnoterule{\vspace*{-5pt}} 
    \begin{center}
    \resizebox{\columnwidth}{!}{%
        \begin{tabular}{ c || l | c | c | c | c | c | c | c | c}
            \toprule
            \textbf{Unit} & \textbf{Cost} & $\boldsymbol{H_1}$ & $\boldsymbol{H_2}$ & $\boldsymbol{H_3}$ & $\boldsymbol{H_4}$ & $\boldsymbol{H_5}$ & $\boldsymbol{H_6}$ & $\boldsymbol{H_7}$ & $\boldsymbol{H_8}$ \\
            \cline{1-10}
            \multirow{4}{*}{$\boldsymbol{G_2}$} & EOM   & $-$ & $-$ & $-$ & $-$ & $-20$ & $-20$ & $-20$ & $-$  \\
                                 & SU    & $-$ & $-$ & $-$ & $-$ & $-300$ & $0$ & $0$ & $-$     \\
                                 & IP    & $-$ & $-$ & $-$ & $-$ & $152$ & $61$ & $144$ & $-$    \\
                                 & Total & $-$ & $-$ & $-$ & $-$ & $-168$ & $41$ & $124$ & $-$         \\
            \cline{1-10}
            \multirow{4}{*}{$\boldsymbol{G_3}$} & EOM   & $-$ & $-$ & $-$ & $-10$ & $-10$ & $-10$ & $-$ & $-$   \\
                                 & SU    & $-$ & $-$ & $-$ & $-200$ & $0$ & $0$ & $-$ & $-$      \\
                                 & IP    & $-$ & $-$ & $-$ & $80$ & $31$ & $122$ & $-$ & $-$  \\
                                 & Total & $-$ & $-$ & $-$ & $-130$ & $21$ & $112$ & $-$ & $-$     \\
            \arrayrulecolor{black}\bottomrule
        \end{tabular}
        }
        \end{center}
    \end{minipage}
    \vspace{-0.35cm}
\end{table}

\subsubsection{Large test system} 
By analyzing the total hourly profits of individual generators given in Table~\ref{tab:SG_profits_M2}, we can conclude that in a larger system more units face negative profits compared to \textit{Method~1}. There are two explanations for such outcome: (i) the issue of symmetry already observed in the small system; (ii) the utility function payment method does not provide compensation for reduced EOM revenue due to a decrease in energy price. Indeed, as shown previously in Fig.~\ref{fig:eom_comp}, the inertia requirements and changed commitment of the system result in lower EOM prices and therefore reduced unit profits. Such effect is not captured by this method and certain generators will inevitably face negative profits and not receive proper reimbursement.

\begin{table}[!b]
\renewcommand{\arraystretch}{0.8}
\vspace{-0.35cm}
\caption{Hourly SG profits [\euro] under \textit{Method~2} pricing scheme.}
\label{tab:SG_profits_M2}
\noindent
\centering
\begin{minipage}{1\linewidth}
\renewcommand\footnoterule{\vspace*{-5pt}}
\begin{center}
    \huge
\resizebox{\columnwidth}{!}{%
\begin{tabular}{ c"c"c"c"c"c"c"c"c"c"c"c"c"c"c"c"c"c"c"c"c}
    \thickhline
    \cellcolor{black!10}$\Huge\boldsymbol{H\text{\textbackslash}G}$& $\Huge\boldsymbol{G_1}$ & $\Huge\boldsymbol{G_2}$ & $\Huge\boldsymbol{G_3}$ & $\Huge\boldsymbol{G_4}$ & $\Huge\boldsymbol{G_5}$ & $\Huge\boldsymbol{G_6}$ & $\Huge\boldsymbol{G_7}$ & $\Huge\boldsymbol{G_8}$ & $\Huge\boldsymbol{G_9}$ & $\Huge\boldsymbol{G_{10}}$ & $\Huge\boldsymbol{G_{11}}$ & $\Huge\boldsymbol{G_{12}}$ & $\Huge\boldsymbol{G_{13}}$ & $\Huge\boldsymbol{G_{14}}$ & $\Huge\boldsymbol{G_{15}}$ & $\Huge\boldsymbol{G_{16}}$ & $\Huge\boldsymbol{G_{17}}$ & $\Huge\boldsymbol{G_{18}}$ & $\Huge\boldsymbol{G_{19}}$ & $\Huge\boldsymbol{G_{20}}$ \\
    \thickhline
$\Huge\boldsymbol{H_1}$ & 	\cellcolor{myGreen!25}168 & 	\cellcolor{myGreen!25}168 & 	\cellcolor{myGreen!25}134 & 	\cellcolor{red!25}-153 & 	\cellcolor{myGreen!25}324 & 	\cellcolor{white}0 & 	\cellcolor{white}0 & 	\cellcolor{white}0 & 	\cellcolor{myGreen!25}364 & 	\cellcolor{red!25}-1222 & 	\cellcolor{myGreen!25}168 & 	\cellcolor{myGreen!25}168 & 	\cellcolor{myGreen!25}134 & 	\cellcolor{red!25}-153 & 	\cellcolor{myGreen!25}324 & 	\cellcolor{white}0 & 	\cellcolor{white}0 & 	\cellcolor{white}0 & 	\cellcolor{myGreen!25}364 & 	\cellcolor{red!25}-784	\\ \thickhline
$\Huge\boldsymbol{H_2}$ & 	\cellcolor{myGreen!25}28 & 	\cellcolor{myGreen!25}28 & 	\cellcolor{myGreen!25}20 & 	\cellcolor{red!25}-321 & 	\cellcolor{myGreen!25}208 & 	\cellcolor{white}0 & 	\cellcolor{white}0 & 	\cellcolor{white}0 & 	\cellcolor{myGreen!25}109 & 	\cellcolor{red!25}-225 & 	\cellcolor{myGreen!25}28 & 	\cellcolor{myGreen!25}28 & 	\cellcolor{myGreen!25}20 & 	\cellcolor{red!25}-321 & 	\cellcolor{myGreen!25}208 & 	\cellcolor{white}0 & 	\cellcolor{white}0 & 	\cellcolor{white}0 & 	\cellcolor{myGreen!25}109 & 	\cellcolor{myGreen!25}79	\\ \thickhline
$\Huge\boldsymbol{H_3}$ & 	\cellcolor{myGreen!25}28 & 	\cellcolor{myGreen!25}28 & 	\cellcolor{myGreen!25}20 & 	\cellcolor{red!25}-321 & 	\cellcolor{myGreen!25}208 & 	\cellcolor{white}0 & 	\cellcolor{white}0 & 	\cellcolor{white}0 & 	\cellcolor{myGreen!25}109 & 	\cellcolor{red!25}-225 & 	\cellcolor{myGreen!25}28 & 	\cellcolor{myGreen!25}28 & 	\cellcolor{myGreen!25}20 & 	\cellcolor{red!25}-321 & 	\cellcolor{myGreen!25}208 & 	\cellcolor{white}0 & 	\cellcolor{white}0 & 	\cellcolor{white}0 & 	\cellcolor{myGreen!25}109 & 	\cellcolor{myGreen!25}79	\\ \thickhline
$\Huge\boldsymbol{H_4}$ & 	\cellcolor{myGreen!25}28 & 	\cellcolor{myGreen!25}28 & 	\cellcolor{myGreen!25}20 & 	\cellcolor{red!25}-321 & 	\cellcolor{myGreen!25}208 & 	\cellcolor{white}0 & 	\cellcolor{white}0 & 	\cellcolor{white}0 & 	\cellcolor{myGreen!25}109 & 	\cellcolor{red!25}-225 & 	\cellcolor{myGreen!25}28 & 	\cellcolor{myGreen!25}28 & 	\cellcolor{myGreen!25}20 & 	\cellcolor{red!25}-321 & 	\cellcolor{myGreen!25}208 & 	\cellcolor{white}0 & 	\cellcolor{white}0 & 	\cellcolor{white}0 & 	\cellcolor{myGreen!25}109 & 	\cellcolor{myGreen!25}79	\\ \thickhline
$\Huge\boldsymbol{H_5}$ & 	\cellcolor{myGreen!25}28 & 	\cellcolor{myGreen!25}28 & 	\cellcolor{myGreen!25}20 & 	\cellcolor{red!25}-321 & 	\cellcolor{myGreen!25}208 & 	\cellcolor{white}0 & 	\cellcolor{white}0 & 	\cellcolor{white}0 & 	\cellcolor{myGreen!25}109 & 	\cellcolor{red!25}-225 & 	\cellcolor{myGreen!25}28 & 	\cellcolor{myGreen!25}28 & 	\cellcolor{myGreen!25}20 & 	\cellcolor{red!25}-321 & 	\cellcolor{myGreen!25}208 & 	\cellcolor{white}0 & 	\cellcolor{white}0 & 	\cellcolor{white}0 & 	\cellcolor{myGreen!25}109 & 	\cellcolor{myGreen!25}79	\\ \thickhline
$\Huge\boldsymbol{H_6}$ & 	\cellcolor{myGreen!25}28 & 	\cellcolor{myGreen!25}28 & 	\cellcolor{myGreen!25}20 & 	\cellcolor{red!25}-321 & 	\cellcolor{myGreen!25}208 & 	\cellcolor{white}0 & 	\cellcolor{white}0 & 	\cellcolor{white}0 & 	\cellcolor{myGreen!25}109 & 	\cellcolor{red!25}-225 & 	\cellcolor{myGreen!25}28 & 	\cellcolor{myGreen!25}28 & 	\cellcolor{myGreen!25}20 & 	\cellcolor{red!25}-321 & 	\cellcolor{myGreen!25}208 & 	\cellcolor{white}0 & 	\cellcolor{white}0 & 	\cellcolor{white}0 & 	\cellcolor{myGreen!25}109 & 	\cellcolor{myGreen!25}79	\\ \thickhline
$\Huge\boldsymbol{H_7}$ & 	\cellcolor{red!25}-52 & 	\cellcolor{red!25}-52 & 	\cellcolor{myGreen!25}79 & 	\cellcolor{red!25}-96 & 	\cellcolor{myGreen!25}20 & 	\cellcolor{white}0 & 	\cellcolor{myGreen!25}3966 & 	\cellcolor{myGreen!25}3966 & 	\cellcolor{red!25}-130 & 	\cellcolor{red!25}-208 & 	\cellcolor{red!25}-52 & 	\cellcolor{red!25}-52 & 	\cellcolor{red!25}-46 & 	\cellcolor{red!25}-161 & 	\cellcolor{myGreen!25}20 & 	\cellcolor{white}0 & 	\cellcolor{myGreen!25}3966 & 	\cellcolor{myGreen!25}3966 & 	\cellcolor{red!25}-130 & 	\cellcolor{red!25}-126	\\ \thickhline
$\Huge\boldsymbol{H_8}$ & 	\cellcolor{red!25}-1 & 	\cellcolor{red!25}-1 & 	\cellcolor{white}0 & 	\cellcolor{myGreen!25}274 & 	\cellcolor{myGreen!25}44 & 	\cellcolor{white}0 & 	\cellcolor{myGreen!25}4424 & 	\cellcolor{myGreen!25}4424 & 	\cellcolor{red!25}-63 & 	\cellcolor{red!25}-107 & 	\cellcolor{red!25}-143 & 	\cellcolor{red!25}-143 & 	\cellcolor{white}0 & 	\cellcolor{myGreen!25}338 & 	\cellcolor{myGreen!25}44 & 	\cellcolor{white}0 & 	\cellcolor{myGreen!25}4424 & 	\cellcolor{myGreen!25}4424 & 	\cellcolor{red!25}-63 & 	\cellcolor{red!25}-8	\\ \thickhline
$\Huge\boldsymbol{H_9}$ & 	\cellcolor{red!25}-106 & 	\cellcolor{red!25}-106 & 	\cellcolor{white}0 & 	\cellcolor{myGreen!25}338 & 	\cellcolor{red!25}-36 & 	\cellcolor{white}0 & 	\cellcolor{myGreen!25}4424 & 	\cellcolor{myGreen!25}4424 & 	\cellcolor{red!25}-228 & 	\cellcolor{red!25}-280 & 	\cellcolor{red!25}-143 & 	\cellcolor{red!25}-143 & 	\cellcolor{white}0 & 	\cellcolor{myGreen!25}338 & 	\cellcolor{red!25}-36 & 	\cellcolor{white}0 & 	\cellcolor{myGreen!25}4424 & 	\cellcolor{myGreen!25}4424 & 	\cellcolor{red!25}-228 & 	\cellcolor{red!25}-254	\\ \thickhline
$\Huge\boldsymbol{H_{10}}$ & 	\cellcolor{red!25}-13 & 	\cellcolor{red!25}-13 & 	\cellcolor{white}0 & 	\cellcolor{myGreen!25}338 & 	\cellcolor{myGreen!25}34 & 	\cellcolor{white}0 & 	\cellcolor{myGreen!25}4424 & 	\cellcolor{myGreen!25}4424 & 	\cellcolor{red!25}-83 & 	\cellcolor{red!25}-128 & 	\cellcolor{red!25}-143 & 	\cellcolor{red!25}-143 & 	\cellcolor{white}0 & 	\cellcolor{myGreen!25}338 & 	\cellcolor{myGreen!25}34 & 	\cellcolor{white}0 & 	\cellcolor{myGreen!25}4424 & 	\cellcolor{myGreen!25}4424 & 	\cellcolor{red!25}-83 & 	\cellcolor{red!25}-38	\\ \thickhline
$\Huge\boldsymbol{H_{11}}$ & 	\cellcolor{myGreen!25}12 & 	\cellcolor{myGreen!25}12 & 	\cellcolor{white}0 & 	\cellcolor{myGreen!25}338 & 	\cellcolor{myGreen!25}53 & 	\cellcolor{white}0 & 	\cellcolor{myGreen!25}4424 & 	\cellcolor{myGreen!25}4424 & 	\cellcolor{red!25}-44 & 	\cellcolor{red!25}-87 & 	\cellcolor{myGreen!25}12 & 	\cellcolor{myGreen!25}12 & 	\cellcolor{white}0 & 	\cellcolor{myGreen!25}338 & 	\cellcolor{myGreen!25}53 & 	\cellcolor{white}0 & 	\cellcolor{myGreen!25}4424 & 	\cellcolor{myGreen!25}4424 & 	\cellcolor{red!25}-44 & 	\cellcolor{myGreen!25}21	\\ \thickhline
$\Huge\boldsymbol{H_{12}}$ & 	\cellcolor{red!25}-1 & 	\cellcolor{red!25}-1 & 	\cellcolor{red!25}-23 & 	\cellcolor{white}0 & 	\cellcolor{myGreen!25}51 & 	\cellcolor{white}0 & 	\cellcolor{myGreen!25}4195 & 	\cellcolor{myGreen!25}4195 & 	\cellcolor{red!25}-57 & 	\cellcolor{red!25}-116 & 	\cellcolor{red!25}-1 & 	\cellcolor{red!25}-1 & 	\cellcolor{myGreen!25}142 & 	\cellcolor{white}0 & 	\cellcolor{myGreen!25}51 & 	\cellcolor{white}0 & 	\cellcolor{myGreen!25}4195 & 	\cellcolor{myGreen!25}4195 & 	\cellcolor{red!25}-57 & 	\cellcolor{red!25}-8	\\ \thickhline
$\Huge\boldsymbol{H_{13}}$ & 	\cellcolor{red!25}-1 & 	\cellcolor{red!25}-1 & 	\cellcolor{red!25}-23 & 	\cellcolor{white}0 & 	\cellcolor{myGreen!25}51 & 	\cellcolor{white}0 & 	\cellcolor{myGreen!25}4195 & 	\cellcolor{myGreen!25}4195 & 	\cellcolor{red!25}-57 & 	\cellcolor{red!25}-116 & 	\cellcolor{red!25}-1 & 	\cellcolor{red!25}-1 & 	\cellcolor{myGreen!25}142 & 	\cellcolor{white}0 & 	\cellcolor{myGreen!25}51 & 	\cellcolor{white}0 & 	\cellcolor{myGreen!25}4195 & 	\cellcolor{myGreen!25}4195 & 	\cellcolor{red!25}-57 & 	\cellcolor{red!25}-8	\\ \thickhline
$\Huge\boldsymbol{H_{14}}$ & 	\cellcolor{red!25}-85 & 	\cellcolor{red!25}-85 & 	\cellcolor{red!25}-442 & 	\cellcolor{red!25}-1091 & 	\cellcolor{myGreen!25}122 & 	\cellcolor{white}0 & 	\cellcolor{white}0 & 	\cellcolor{white}0 & 	\cellcolor{red!25}-67 & 	\cellcolor{red!25}-410 & 	\cellcolor{red!25}-85 & 	\cellcolor{red!25}-85 & 	\cellcolor{red!25}-100 & 	\cellcolor{red!25}-1091 & 	\cellcolor{myGreen!25}122 & 	\cellcolor{white}0 & 	\cellcolor{white}0 & 	\cellcolor{white}0 & 	\cellcolor{red!25}-67 & 	\cellcolor{red!25}-185	\\ \thickhline
$\Huge\boldsymbol{H_{15}}$ & 	\cellcolor{myGreen!25}30 & 	\cellcolor{myGreen!25}30 & 	\cellcolor{red!25}-442 & 	\cellcolor{red!25}-319 & 	\cellcolor{myGreen!25}209 & 	\cellcolor{white}0 & 	\cellcolor{white}0 & 	\cellcolor{white}0 & 	\cellcolor{myGreen!25}111 & 	\cellcolor{red!25}-223 & 	\cellcolor{myGreen!25}30 & 	\cellcolor{myGreen!25}30 & 	\cellcolor{myGreen!25}22 & 	\cellcolor{red!25}-319 & 	\cellcolor{myGreen!25}209 & 	\cellcolor{white}0 & 	\cellcolor{white}0 & 	\cellcolor{white}0 & 	\cellcolor{myGreen!25}111 & 	\cellcolor{myGreen!25}82	\\ \thickhline
$\Huge\boldsymbol{H_{16}}$ & 	\cellcolor{red!25}-23 & 	\cellcolor{red!25}-23 & 	\cellcolor{red!25}-442 & 	\cellcolor{red!25}-412 & 	\cellcolor{myGreen!25}169 & 	\cellcolor{white}0 & 	\cellcolor{white}0 & 	\cellcolor{white}0 & 	\cellcolor{myGreen!25}29 & 	\cellcolor{red!25}-309 & 	\cellcolor{red!25}-23 & 	\cellcolor{red!25}-23 & 	\cellcolor{red!25}-34 & 	\cellcolor{red!25}-412 & 	\cellcolor{myGreen!25}169 & 	\cellcolor{white}0 & 	\cellcolor{white}0 & 	\cellcolor{white}0 & 	\cellcolor{myGreen!25}29 & 	\cellcolor{red!25}-41	\\ \thickhline
$\Huge\boldsymbol{H_{17}}$ & 	\cellcolor{red!25}-29 & 	\cellcolor{red!25}-29 & 	\cellcolor{red!25}-442 & 	\cellcolor{red!25}-422 & 	\cellcolor{myGreen!25}164 & 	\cellcolor{white}0 & 	\cellcolor{white}0 & 	\cellcolor{white}0 & 	\cellcolor{myGreen!25}20 & 	\cellcolor{red!25}-319 & 	\cellcolor{red!25}-29 & 	\cellcolor{red!25}-29 & 	\cellcolor{red!25}-40 & 	\cellcolor{red!25}-422 & 	\cellcolor{myGreen!25}164 & 	\cellcolor{white}0 & 	\cellcolor{white}0 & 	\cellcolor{white}0 & 	\cellcolor{myGreen!25}20 & 	\cellcolor{red!25}-55	\\ \thickhline
$\Huge\boldsymbol{H_{18}}$ & 	\cellcolor{myGreen!25}28 & 	\cellcolor{myGreen!25}28 & 	\cellcolor{myGreen!25}20 & 	\cellcolor{red!25}-321 & 	\cellcolor{myGreen!25}208 & 	\cellcolor{white}0 & 	\cellcolor{white}0 & 	\cellcolor{white}0 & 	\cellcolor{myGreen!25}109 & 	\cellcolor{red!25}-225 & 	\cellcolor{myGreen!25}28 & 	\cellcolor{myGreen!25}28 & 	\cellcolor{myGreen!25}20 & 	\cellcolor{red!25}-321 & 	\cellcolor{myGreen!25}208 & 	\cellcolor{white}0 & 	\cellcolor{white}0 & 	\cellcolor{white}0 & 	\cellcolor{myGreen!25}109 & 	\cellcolor{myGreen!25}79	\\ \thickhline
$\Huge\boldsymbol{H_{19}}$ & 	\cellcolor{myGreen!25}28 & 	\cellcolor{myGreen!25}28 & 	\cellcolor{myGreen!25}20 & 	\cellcolor{red!25}-321 & 	\cellcolor{myGreen!25}208 & 	\cellcolor{white}0 & 	\cellcolor{white}0 & 	\cellcolor{white}0 & 	\cellcolor{myGreen!25}109 & 	\cellcolor{red!25}-225 & 	\cellcolor{myGreen!25}28 & 	\cellcolor{myGreen!25}28 & 	\cellcolor{myGreen!25}20 & 	\cellcolor{red!25}-321 & 	\cellcolor{myGreen!25}208 & 	\cellcolor{white}0 & 	\cellcolor{white}0 & 	\cellcolor{white}0 & 	\cellcolor{myGreen!25}109 & 	\cellcolor{myGreen!25}79	\\ \thickhline
$\Huge\boldsymbol{H_{20}}$ & 	\cellcolor{myGreen!25}8 & 	\cellcolor{myGreen!25}8 & 	\cellcolor{red!25}-107 & 	\cellcolor{red!25}-804 & 	\cellcolor{myGreen!25}285 & 	\cellcolor{white}0 & 	\cellcolor{red!25}-723 & 	\cellcolor{red!25}-723 & 	\cellcolor{myGreen!25}160 & 	\cellcolor{red!25}-367 & 	\cellcolor{myGreen!25}8 & 	\cellcolor{myGreen!25}8 & 	\cellcolor{red!25}-107 & 	\cellcolor{red!25}-804 & 	\cellcolor{myGreen!25}285 & 	\cellcolor{white}0 & 	\cellcolor{red!25}-723 & 	\cellcolor{red!25}-723 & 	\cellcolor{myGreen!25}160 & 	\cellcolor{myGreen!25}44	\\ \thickhline
$\Huge\boldsymbol{H_{21}}$ & 	\cellcolor{myGreen!25}75 & 	\cellcolor{myGreen!25}75 & 	\cellcolor{red!25}-36 & 	\cellcolor{red!25}-686 & 	\cellcolor{myGreen!25}336 & 	\cellcolor{white}0 & 	\cellcolor{red!25}-723 & 	\cellcolor{red!25}-723 & 	\cellcolor{myGreen!25}265 & 	\cellcolor{red!25}-258 & 	\cellcolor{myGreen!25}75 & 	\cellcolor{myGreen!25}75 & 	\cellcolor{red!25}-36 & 	\cellcolor{red!25}-686 & 	\cellcolor{myGreen!25}336 & 	\cellcolor{white}0 & 	\cellcolor{red!25}-723 & 	\cellcolor{red!25}-723 & 	\cellcolor{myGreen!25}265 & 	\cellcolor{myGreen!25}200	\\ \thickhline
$\Huge\boldsymbol{H_{22}}$ & 	\cellcolor{myGreen!25}75 & 	\cellcolor{myGreen!25}75 & 	\cellcolor{red!25}-36 & 	\cellcolor{red!25}-686 & 	\cellcolor{myGreen!25}336 & 	\cellcolor{white}0 & 	\cellcolor{red!25}-723 & 	\cellcolor{red!25}-723 & 	\cellcolor{myGreen!25}265 & 	\cellcolor{red!25}-258 & 	\cellcolor{myGreen!25}75 & 	\cellcolor{myGreen!25}75 & 	\cellcolor{red!25}-36 & 	\cellcolor{red!25}-686 & 	\cellcolor{myGreen!25}336 & 	\cellcolor{white}0 & 	\cellcolor{red!25}-723 & 	\cellcolor{red!25}-723 & 	\cellcolor{myGreen!25}265 & 	\cellcolor{myGreen!25}200	\\ \thickhline
$\Huge\boldsymbol{H_{23}}$ & 	\cellcolor{myGreen!25}75 & 	\cellcolor{myGreen!25}75 & 	\cellcolor{red!25}-36 & 	\cellcolor{red!25}-686 & 	\cellcolor{myGreen!25}336 & 	\cellcolor{white}0 & 	\cellcolor{red!25}-723 & 	\cellcolor{red!25}-723 & 	\cellcolor{myGreen!25}265 & 	\cellcolor{red!25}-258 & 	\cellcolor{myGreen!25}75 & 	\cellcolor{myGreen!25}75 & 	\cellcolor{red!25}-36 & 	\cellcolor{red!25}-686 & 	\cellcolor{myGreen!25}336 & 	\cellcolor{white}0 & 	\cellcolor{red!25}-723 & 	\cellcolor{red!25}-723 & 	\cellcolor{myGreen!25}265 & 	\cellcolor{myGreen!25}200	\\ \thickhline
$\Huge\boldsymbol{H_{24}}$ & 	\cellcolor{myGreen!25}75 & 	\cellcolor{myGreen!25}75 & 	\cellcolor{red!25}-36 & 	\cellcolor{red!25}-686 & 	\cellcolor{myGreen!25}336 & 	\cellcolor{white}0 & 	\cellcolor{red!25}-723 & 	\cellcolor{red!25}-723 & 	\cellcolor{myGreen!25}265 & 	\cellcolor{red!25}-258 & 	\cellcolor{myGreen!25}75 & 	\cellcolor{myGreen!25}75 & 	\cellcolor{red!25}-36 & 	\cellcolor{red!25}-686 & 	\cellcolor{myGreen!25}336 & 	\cellcolor{white}0 & 	\cellcolor{red!25}-723 & 	\cellcolor{red!25}-723 & 	\cellcolor{myGreen!25}265 & 	\cellcolor{myGreen!25}200	\\ 
    \arrayrulecolor{black}\thickhline
\end{tabular}
 }
\end{center}
\end{minipage}
\end{table}

The reason why payments to inertia providers lead to some generators making the exact opposite profits lies in the fact that the TSO reimburses them according to the exact amount at which it values the inertia at each hour, which is the cost of having those generators online in case no cheaper inertia is available. There is thus a single price for inertia at each hour, which does not necessarily correspond to the marginal inertia price (i.e., to the most expensive unit online), but to a total budget that will be distributed among providers. Considering that inertia is provided in discrete steps due to the binary nature of the UC problem, the total available inertia of additionally committed generators might end up being higher than the inertia demand $H_\mathrm{dem}$ used in~\eqref{eq:utility}. In such instances, the allocation of inertia to each generator is arbitrary and thus some units might receive disproportionate payments compared to others during the redistribution of inertia provision budget. This problem can be resolved by assigning the exact amount of inertia provided by each generator \textit{ex-post}, such that no generator is facing losses. Nevertheless, it is hard to find a systematic way in which this procedure should be conducted. This is in fact a main distinction compared to \textit{Method~1}, where the price of inertia is calculated \textit{ex-post} as the cost of the most expensive unit coming online for inertia purposes.

\subsubsection{Inclusion of Virtual Inertia}
The addition of virtual inertia creates issues for reimbursement under the utility function method. More precisely, with the addition of VI units the dual of the RoCoF constraint yields $\lambda_t^\mathrm{H} = U_t^\mathrm{H} - C_v^\mathrm{VI}$, which could have negative values for large enough $C_v^\mathrm{VI}$ and result in units paying to provide inertia. In other words, as virtual inertia becomes more expensive, less payments will be made to units offering inertia which goes against the basic market intuition. Such phenomena can also be seen in Fig.~\ref{fig:VI_payments} on the example of a small test system and assuming all VI units bid-in at the same cost. Furthermore, for $C_v^\mathrm{VI}\in\left[150,200\right]\,\text{\euro}/\mathrm{MWs}^2$ the dual becomes negative which is conceptually problematic.

A possible solution could simply be to set the dual to $\lambda_t^\mathrm{H} = U_t^\mathrm{H}$, with the TSO always paying $U^\mathrm{H}$ regardless of which technology is providing inertia and at which price. This would encourage cheaper technologies to enter the market with the potential for high profits. Such solution is not necessarily a disadvantage for the TSO. As a matter of fact, the operator could opt for paying a fixed price for inertia in order to encourage units with cheaper inertia provision to penetrate the system; once this is achieved, the TSO can reduce inertia payments.

\begin{table}[!t]
\renewcommand{\arraystretch}{0.8}
\caption{Hourly SG profits [\euro] under \textit{Method~3} pricing scheme.}
\label{tab:SG_profits_M3}
\noindent
\centering
\begin{minipage}{0.96\linewidth}
\renewcommand\footnoterule{\vspace*{-5pt}}
\begin{center}
    \huge
\resizebox{\columnwidth}{!}{%
\begin{tabular}{ c"c"c"c"c"c"c"c"c"c"c"c"c"c"c"c"c"c"c"c"c}
    \thickhline
    \cellcolor{black!10}$\Huge\boldsymbol{H\text{\textbackslash}G}$& $\Huge\boldsymbol{G_1}$ & $\Huge\boldsymbol{G_2}$ & $\Huge\boldsymbol{G_3}$ & $\Huge\boldsymbol{G_4}$ & $\Huge\boldsymbol{G_5}$ & $\Huge\boldsymbol{G_6}$ & $\Huge\boldsymbol{G_7}$ & $\Huge\boldsymbol{G_8}$ & $\Huge\boldsymbol{G_9}$ & $\Huge\boldsymbol{G_{10}}$ & $\Huge\boldsymbol{G_{11}}$ & $\Huge\boldsymbol{G_{12}}$ & $\Huge\boldsymbol{G_{13}}$ & $\Huge\boldsymbol{G_{14}}$ & $\Huge\boldsymbol{G_{15}}$ & $\Huge\boldsymbol{G_{16}}$ & $\Huge\boldsymbol{G_{17}}$ & $\Huge\boldsymbol{G_{18}}$ & $\Huge\boldsymbol{G_{19}}$ & $\Huge\boldsymbol{G_{20}}$ \\
    \thickhline
$\Huge\boldsymbol{H_1}$ & 	\cellcolor{white}0 & 	\cellcolor{white}0 & 	\cellcolor{white}0 & 	\cellcolor{white}0 & 	\cellcolor{white}0 & 	\cellcolor{white}0 & 	\cellcolor{white}0 & 	\cellcolor{white}0 & 	\cellcolor{white}0 & 	\cellcolor{white}0 & 	\cellcolor{white}0 & 	\cellcolor{white}0 & 	\cellcolor{white}0 & 	\cellcolor{white}0 & 	\cellcolor{white}0 & 	\cellcolor{white}0 & 	\cellcolor{white}0 & 	\cellcolor{white}0 & 	\cellcolor{white}0 & 	\cellcolor{white}0	\\ \thickhline
$\Huge\boldsymbol{H_2}$ & 	\cellcolor{white}0 & 	\cellcolor{white}0 & 	\cellcolor{white}0 & 	\cellcolor{white}0 & 	\cellcolor{white}0 & 	\cellcolor{white}0 & 	\cellcolor{white}0 & 	\cellcolor{white}0 & 	\cellcolor{white}0 & 	\cellcolor{white}0 & 	\cellcolor{white}0 & 	\cellcolor{white}0 & 	\cellcolor{white}0 & 	\cellcolor{white}0 & 	\cellcolor{white}0 & 	\cellcolor{white}0 & 	\cellcolor{white}0 & 	\cellcolor{white}0 & 	\cellcolor{white}0 & 	\cellcolor{white}0	\\ \thickhline
$\Huge\boldsymbol{H_3}$ & 	\cellcolor{white}0 & 	\cellcolor{white}0 & 	\cellcolor{white}0 & 	\cellcolor{white}0 & 	\cellcolor{white}0 & 	\cellcolor{white}0 & 	\cellcolor{white}0 & 	\cellcolor{white}0 & 	\cellcolor{white}0 & 	\cellcolor{white}0 & 	\cellcolor{white}0 & 	\cellcolor{white}0 & 	\cellcolor{white}0 & 	\cellcolor{white}0 & 	\cellcolor{white}0 & 	\cellcolor{white}0 & 	\cellcolor{white}0 & 	\cellcolor{white}0 & 	\cellcolor{white}0 & 	\cellcolor{white}0	\\ \thickhline
$\Huge\boldsymbol{H_4}$ & 	\cellcolor{white}0 & 	\cellcolor{white}0 & 	\cellcolor{white}0 & 	\cellcolor{white}0 & 	\cellcolor{white}0 & 	\cellcolor{white}0 & 	\cellcolor{white}0 & 	\cellcolor{white}0 & 	\cellcolor{white}0 & 	\cellcolor{white}0 & 	\cellcolor{white}0 & 	\cellcolor{white}0 & 	\cellcolor{white}0 & 	\cellcolor{white}0 & 	\cellcolor{white}0 & 	\cellcolor{white}0 & 	\cellcolor{white}0 & 	\cellcolor{white}0 & 	\cellcolor{white}0 & 	\cellcolor{white}0	\\ \thickhline
$\Huge\boldsymbol{H_5}$ & 	\cellcolor{white}0 & 	\cellcolor{white}0 & 	\cellcolor{white}0 & 	\cellcolor{white}0 & 	\cellcolor{white}0 & 	\cellcolor{white}0 & 	\cellcolor{white}0 & 	\cellcolor{white}0 & 	\cellcolor{white}0 & 	\cellcolor{white}0 & 	\cellcolor{white}0 & 	\cellcolor{white}0 & 	\cellcolor{white}0 & 	\cellcolor{white}0 & 	\cellcolor{white}0 & 	\cellcolor{white}0 & 	\cellcolor{white}0 & 	\cellcolor{white}0 & 	\cellcolor{white}0 & 	\cellcolor{white}0	\\ \thickhline
$\Huge\boldsymbol{H_6}$ & 	\cellcolor{white}0 & 	\cellcolor{white}0 & 	\cellcolor{white}0 & 	\cellcolor{white}0 & 	\cellcolor{white}0 & 	\cellcolor{white}0 & 	\cellcolor{white}0 & 	\cellcolor{white}0 & 	\cellcolor{white}0 & 	\cellcolor{white}0 & 	\cellcolor{white}0 & 	\cellcolor{white}0 & 	\cellcolor{white}0 & 	\cellcolor{white}0 & 	\cellcolor{white}0 & 	\cellcolor{white}0 & 	\cellcolor{white}0 & 	\cellcolor{white}0 & 	\cellcolor{white}0 & 	\cellcolor{white}0	\\ \thickhline
$\Huge\boldsymbol{H_7}$ & 	\cellcolor{white}0 & 	\cellcolor{white}0 & 	\cellcolor{white}0 & 	\cellcolor{white}0 & 	\cellcolor{white}0 & 	\cellcolor{white}0 & 	\cellcolor{myGreen!25}3966 & 	\cellcolor{myGreen!25}3966 & 	\cellcolor{white}0 & 	\cellcolor{white}0 & 	\cellcolor{white}0 & 	\cellcolor{white}0 & 	\cellcolor{white}0 & 	\cellcolor{white}0 & 	\cellcolor{white}0 & 	\cellcolor{white}0 & 	\cellcolor{myGreen!25}3966 & 	\cellcolor{myGreen!25}3966 & 	\cellcolor{white}0 & 	\cellcolor{white}0	\\ \thickhline
$\Huge\boldsymbol{H_8}$ & 	\cellcolor{white}0 & 	\cellcolor{white}0 & 	\cellcolor{white}0 & 	\cellcolor{myGreen!25}338 & 	\cellcolor{white}0 & 	\cellcolor{white}0 & 	\cellcolor{myGreen!25}4424 & 	\cellcolor{myGreen!25}4424 & 	\cellcolor{white}0 & 	\cellcolor{white}0 & 	\cellcolor{white}0 & 	\cellcolor{white}0 & 	\cellcolor{white}0 & 	\cellcolor{myGreen!25}274 & 	\cellcolor{white}0 & 	\cellcolor{white}0 & 	\cellcolor{myGreen!25}4424 & 	\cellcolor{myGreen!25}4424 & 	\cellcolor{white}0 & 	\cellcolor{white}0	\\ \thickhline
$\Huge\boldsymbol{H_9}$ & 	\cellcolor{white}0 & 	\cellcolor{white}0 & 	\cellcolor{white}0 & 	\cellcolor{myGreen!25}338 & 	\cellcolor{white}0 & 	\cellcolor{white}0 & 	\cellcolor{myGreen!25}4424 & 	\cellcolor{myGreen!25}4424 & 	\cellcolor{white}0 & 	\cellcolor{white}0 & 	\cellcolor{white}0 & 	\cellcolor{white}0 & 	\cellcolor{white}0 & 	\cellcolor{myGreen!25}338 & 	\cellcolor{white}0 & 	\cellcolor{white}0 & 	\cellcolor{myGreen!25}4424 & 	\cellcolor{myGreen!25}4424 & 	\cellcolor{white}0 & 	\cellcolor{white}0	\\ \thickhline
$\Huge\boldsymbol{H_{10}}$ & 	\cellcolor{white}0 & 	\cellcolor{white}0 & 	\cellcolor{white}0 & 	\cellcolor{myGreen!25}338 & 	\cellcolor{white}0 & 	\cellcolor{white}0 & 	\cellcolor{myGreen!25}4424 & 	\cellcolor{myGreen!25}4424 & 	\cellcolor{white}0 & 	\cellcolor{white}0 & 	\cellcolor{white}0 & 	\cellcolor{white}0 & 	\cellcolor{white}0 & 	\cellcolor{myGreen!25}338 & 	\cellcolor{white}0 & 	\cellcolor{white}0 & 	\cellcolor{myGreen!25}4424 & 	\cellcolor{myGreen!25}4424 & 	\cellcolor{white}0 & 	\cellcolor{white}0	\\ \thickhline
$\Huge\boldsymbol{H_{11}}$ & 	\cellcolor{white}0 & 	\cellcolor{white}0 & 	\cellcolor{white}0 & 	\cellcolor{myGreen!25}338 & 	\cellcolor{white}0 & 	\cellcolor{white}0 & 	\cellcolor{myGreen!25}4424 & 	\cellcolor{myGreen!25}4424 & 	\cellcolor{white}0 & 	\cellcolor{white}0 & 	\cellcolor{white}0 & 	\cellcolor{white}0 & 	\cellcolor{white}0 & 	\cellcolor{myGreen!25}338 & 	\cellcolor{white}0 & 	\cellcolor{white}0 & 	\cellcolor{myGreen!25}4424 & 	\cellcolor{myGreen!25}4424 & 	\cellcolor{white}0 & 	\cellcolor{white}0	\\ \thickhline
$\Huge\boldsymbol{H_{12}}$ & 	\cellcolor{white}0 & 	\cellcolor{white}0 & 	\cellcolor{white}0 & 	\cellcolor{white}0 & 	\cellcolor{white}0 & 	\cellcolor{white}0 & 	\cellcolor{myGreen!25}4195 & 	\cellcolor{myGreen!25}4195 & 	\cellcolor{white}0 & 	\cellcolor{white}0 & 	\cellcolor{white}0 & 	\cellcolor{white}0 & 	\cellcolor{white}0 & 	\cellcolor{white}0 & 	\cellcolor{white}0 & 	\cellcolor{white}0 & 	\cellcolor{myGreen!25}4195 & 	\cellcolor{myGreen!25}4195 & 	\cellcolor{white}0 & 	\cellcolor{white}0	\\ \thickhline
$\Huge\boldsymbol{H_{13}}$ & 	\cellcolor{white}0 & 	\cellcolor{white}0 & 	\cellcolor{white}0 & 	\cellcolor{white}0 & 	\cellcolor{white}0 & 	\cellcolor{white}0 & 	\cellcolor{myGreen!25}4195 & 	\cellcolor{myGreen!25}4195 & 	\cellcolor{white}0 & 	\cellcolor{white}0 & 	\cellcolor{white}0 & 	\cellcolor{white}0 & 	\cellcolor{white}0 & 	\cellcolor{white}0 & 	\cellcolor{white}0 & 	\cellcolor{white}0 & 	\cellcolor{myGreen!25}4195 & 	\cellcolor{myGreen!25}4195 & 	\cellcolor{white}0 & 	\cellcolor{white}0	\\ \thickhline
$\Huge\boldsymbol{H_{14}}$ & 	\cellcolor{white}0 & 	\cellcolor{white}0 & 	\cellcolor{white}0 & 	\cellcolor{white}0 & 	\cellcolor{white}0 & 	\cellcolor{white}0 & 	\cellcolor{white}0 & 	\cellcolor{white}0 & 	\cellcolor{white}0 & 	\cellcolor{white}0 & 	\cellcolor{white}0 & 	\cellcolor{white}0 & 	\cellcolor{white}0 & 	\cellcolor{white}0 & 	\cellcolor{white}0 & 	\cellcolor{white}0 & 	\cellcolor{white}0 & 	\cellcolor{white}0 & 	\cellcolor{white}0 & 	\cellcolor{white}0	\\ \thickhline
$\Huge\boldsymbol{H_{15}}$ & 	\cellcolor{white}0 & 	\cellcolor{white}0 & 	\cellcolor{white}0 & 	\cellcolor{white}0 & 	\cellcolor{white}0 & 	\cellcolor{white}0 & 	\cellcolor{white}0 & 	\cellcolor{white}0 & 	\cellcolor{white}0 & 	\cellcolor{white}0 & 	\cellcolor{white}0 & 	\cellcolor{white}0 & 	\cellcolor{white}0 & 	\cellcolor{white}0 & 	\cellcolor{white}0 & 	\cellcolor{white}0 & 	\cellcolor{white}0 & 	\cellcolor{white}0 & 	\cellcolor{white}0 & 	\cellcolor{white}0	\\ \thickhline
$\Huge\boldsymbol{H_{16}}$ & 	\cellcolor{white}0 & 	\cellcolor{white}0 & 	\cellcolor{white}0 & 	\cellcolor{white}0 & 	\cellcolor{white}0 & 	\cellcolor{white}0 & 	\cellcolor{white}0 & 	\cellcolor{white}0 & 	\cellcolor{white}0 & 	\cellcolor{white}0 & 	\cellcolor{white}0 & 	\cellcolor{white}0 & 	\cellcolor{white}0 & 	\cellcolor{white}0 & 	\cellcolor{white}0 & 	\cellcolor{white}0 & 	\cellcolor{white}0 & 	\cellcolor{white}0 & 	\cellcolor{white}0 & 	\cellcolor{white}0	\\ \thickhline
$\Huge\boldsymbol{H_{17}}$ & 	\cellcolor{white}0 & 	\cellcolor{white}0 & 	\cellcolor{white}0 & 	\cellcolor{white}0 & 	\cellcolor{white}0 & 	\cellcolor{white}0 & 	\cellcolor{white}0 & 	\cellcolor{white}0 & 	\cellcolor{white}0 & 	\cellcolor{white}0 & 	\cellcolor{white}0 & 	\cellcolor{white}0 & 	\cellcolor{white}0 & 	\cellcolor{white}0 & 	\cellcolor{white}0 & 	\cellcolor{white}0 & 	\cellcolor{white}0 & 	\cellcolor{white}0 & 	\cellcolor{white}0 & 	\cellcolor{white}0	\\ \thickhline
$\Huge\boldsymbol{H_{18}}$ & 	\cellcolor{white}0 & 	\cellcolor{white}0 & 	\cellcolor{white}0 & 	\cellcolor{white}0 & 	\cellcolor{white}0 & 	\cellcolor{white}0 & 	\cellcolor{white}0 & 	\cellcolor{white}0 & 	\cellcolor{white}0 & 	\cellcolor{white}0 & 	\cellcolor{white}0 & 	\cellcolor{white}0 & 	\cellcolor{white}0 & 	\cellcolor{white}0 & 	\cellcolor{white}0 & 	\cellcolor{white}0 & 	\cellcolor{white}0 & 	\cellcolor{white}0 & 	\cellcolor{white}0 & 	\cellcolor{white}0	\\ \thickhline
$\Huge\boldsymbol{H_{19}}$ & 	\cellcolor{white}0 & 	\cellcolor{white}0 & 	\cellcolor{white}0 & 	\cellcolor{white}0 & 	\cellcolor{white}0 & 	\cellcolor{white}0 & 	\cellcolor{white}0 & 	\cellcolor{white}0 & 	\cellcolor{white}0 & 	\cellcolor{white}0 & 	\cellcolor{white}0 & 	\cellcolor{white}0 & 	\cellcolor{white}0 & 	\cellcolor{white}0 & 	\cellcolor{white}0 & 	\cellcolor{white}0 & 	\cellcolor{white}0 & 	\cellcolor{white}0 & 	\cellcolor{white}0 & 	\cellcolor{white}0	\\ \thickhline
$\Huge\boldsymbol{H_{20}}$ & 	\cellcolor{white}0 & 	\cellcolor{white}0 & 	\cellcolor{white}0 & 	\cellcolor{white}0 & 	\cellcolor{white}0 & 	\cellcolor{white}0 & 	\cellcolor{white}0 & 	\cellcolor{white}0 & 	\cellcolor{white}0 & 	\cellcolor{white}0 & 	\cellcolor{white}0 & 	\cellcolor{white}0 & 	\cellcolor{white}0 & 	\cellcolor{white}0 & 	\cellcolor{white}0 & 	\cellcolor{white}0 & 	\cellcolor{white}0 & 	\cellcolor{white}0 & 	\cellcolor{white}0 & 	\cellcolor{white}0	\\ \thickhline
$\Huge\boldsymbol{H_{21}}$ & 	\cellcolor{white}0 & 	\cellcolor{white}0 & 	\cellcolor{white}0 & 	\cellcolor{white}0 & 	\cellcolor{white}0 & 	\cellcolor{white}0 & 	\cellcolor{white}0 & 	\cellcolor{white}0 & 	\cellcolor{white}0 & 	\cellcolor{white}0 & 	\cellcolor{white}0 & 	\cellcolor{white}0 & 	\cellcolor{white}0 & 	\cellcolor{white}0 & 	\cellcolor{white}0 & 	\cellcolor{white}0 & 	\cellcolor{white}0 & 	\cellcolor{white}0 & 	\cellcolor{white}0 & 	\cellcolor{white}0	\\ \thickhline
$\Huge\boldsymbol{H_{22}}$ & 	\cellcolor{white}0 & 	\cellcolor{white}0 & 	\cellcolor{white}0 & 	\cellcolor{white}0 & 	\cellcolor{white}0 & 	\cellcolor{white}0 & 	\cellcolor{white}0 & 	\cellcolor{white}0 & 	\cellcolor{white}0 & 	\cellcolor{white}0 & 	\cellcolor{white}0 & 	\cellcolor{white}0 & 	\cellcolor{white}0 & 	\cellcolor{white}0 & 	\cellcolor{white}0 & 	\cellcolor{white}0 & 	\cellcolor{white}0 & 	\cellcolor{white}0 & 	\cellcolor{white}0 & 	\cellcolor{white}0	\\ \thickhline
$\Huge\boldsymbol{H_{23}}$ & 	\cellcolor{white}0 & 	\cellcolor{white}0 & 	\cellcolor{white}0 & 	\cellcolor{white}0 & 	\cellcolor{white}0 & 	\cellcolor{white}0 & 	\cellcolor{white}0 & 	\cellcolor{white}0 & 	\cellcolor{white}0 & 	\cellcolor{white}0 & 	\cellcolor{white}0 & 	\cellcolor{white}0 & 	\cellcolor{white}0 & 	\cellcolor{white}0 & 	\cellcolor{white}0 & 	\cellcolor{white}0 & 	\cellcolor{white}0 & 	\cellcolor{white}0 & 	\cellcolor{white}0 & 	\cellcolor{white}0	\\ \thickhline
$\Huge\boldsymbol{H_{24}}$ & 	\cellcolor{white}0 & 	\cellcolor{white}0 & 	\cellcolor{white}0 & 	\cellcolor{white}0 & 	\cellcolor{white}0 & 	\cellcolor{white}0 & 	\cellcolor{white}0 & 	\cellcolor{white}0 & 	\cellcolor{white}0 & 	\cellcolor{white}0 & 	\cellcolor{white}0 & 	\cellcolor{white}0 & 	\cellcolor{white}0 & 	\cellcolor{white}0 & 	\cellcolor{white}0 & 	\cellcolor{white}0 & 	\cellcolor{white}0 & 	\cellcolor{white}0 & 	\cellcolor{white}0 & 	\cellcolor{white}0	\\
    \arrayrulecolor{black}\thickhline
\end{tabular}
 }
\end{center}
\end{minipage}
\vspace{-0.35cm}
\end{table}

\setcounter{table}{10}
\begin{table*}[]
\renewcommand{\arraystretch}{1.2}
\caption{Policy implications of proposed payment methods with distinctive advantages (green) and disadvantages (red).}
\label{tab:methods_summary}
\noindent
\centering
    \begin{minipage}{\linewidth} 
    \renewcommand\footnoterule{\vspace*{-5pt}} 
    \begin{center}
        \begin{tabular}{ c || >{\centering\arraybackslash}m{3.5cm} | >{\centering\arraybackslash}m{6.5cm} | >{\centering\arraybackslash}m{4.5cm} }
            \toprule
            \textbf{Perspective} & \textbf{Method 1: Ex-Post Price} & \textbf{Method 2: Utility Function} & \textbf{Method 3: Uplift Payments}\\
            \cline{1-4}
            \multirow{4.15}{*}{\textbf{TSO}} & {\color{myGreen}Incentives for provision of virtual inertia} & {\color{myGreen}Clear financial expectations due to fixed rates for inertia services} & {\color{myGreen}Low payments whenever more affordable inertia is available} \\
            & {\color{myRed}Large inertia payments} & {\color{myRed}Same expenses regardless of the available inertia resources, unless tailor-made utility function is used} & {\color{myRed}No incentives for investment in more efficient inertia technologies}  \\
            \cline{1-4}
            \multirow{2.25}{*}{\textbf{SGs}} & {\color{myGreen} High profits} & {\color{myGreen}$-$} & {\color{myGreen}$-$} \\
            & {\color{myRed}$-$} & {\color{myRed}Possible negative profits} & {\color{myRed}$-$}  \\
            \cline{1-4}
            \multirow{2.25}{*}{\textbf{VI units}} & {\color{myGreen} High profits} & {\color{myGreen} Clear income expectations} & {\color{myGreen} Clear market understanding} \\
            & {\color{myRed}$-$} & {\color{myRed}Fixed income for inertia provision} & {\color{myRed}Low profitability}  \\
            \arrayrulecolor{black}\bottomrule
        \end{tabular}
        \end{center}
    \end{minipage}
    \vspace{-0.35cm}
\end{table*}

\vspace{-0.35cm}
\subsection{Method 3: Uplift Payments}

Applying the uplift method to the small system yields total inertia payments of \SI{590}{\sieuro}, with losses on the EOM and start-up costs both covered. Therefore, $G_2$ and $G_3$ are left with zero profits. Similar conclusions can be drawn from the hourly profit analysis for the larger test case given in Table~\ref{tab:SG_profits_M3}, where negative profits of all SGs have been reduced to zero.

Adding virtual inertia does not have any negative impact on the payments through \textit{Method~3}. SGs are still compensated to ensure cost recovery and batteries are reimbursed for VI provision based on their bids. Inertia payments will always be lower compared to the case without virtual inertia, as indicated in Fig.~\ref{fig:VI_payments}.

\vspace{-0.35cm}
\subsection{Comparison of Payment Methods} 

\setcounter{table}{8}
\begin{table}[!t]
\renewcommand{\arraystretch}{1.2}
\caption{Comparison of different inertia payment methods.}
\label{tab:comp_methods}
\noindent
\centering
    \begin{minipage}{\linewidth} 
    \renewcommand\footnoterule{\vspace*{-5pt}} 
    \begin{center}
        \begin{tabular}{ l || c | c | c }
            \toprule
            \textbf{System performance metric} & \textbf{Method 1} & \textbf{Method 2} & \textbf{Method 3} \\ 
            \cline{1-4}
            Optimal UC solution [\euro] & $569\,610$ & $569\,610$ & $569\,610$ \\
            Total inertia payments [\euro] & $258\,160$ & $162\,124$ & $187\,948$ \\
            \# units with negative profits & $0/18$  & $5/18$ & $0/18$ \\
            \# units with positive profits & $18/18$ & $13/18$ & $6/18$ \\
            \arrayrulecolor{black}\bottomrule
        \end{tabular}
        \end{center}
    \end{minipage}
    \vspace{-0.35cm}
\end{table}

We now compare the three payment methods and evaluate them on a large test system based on several metrics of performance. Table~\ref{tab:comp_methods} shows the comparison of total system costs and inertia payments for each pricing scheme, as well as the number of units facing negative or positive profits. Understandably, all methods result in the same UC costs and distinguish only between the amount of final inertia payments. Moreover, they yield the same commitment and dispatch schedule for all generators. \textit{Method~2} is the most affordable for the TSO, with lower inertia payments compared to the other two mechanisms. 

\setcounter{table}{9}
\begin{table}[!b]
\renewcommand{\arraystretch}{1.2}
\vspace{-0.35cm}
\caption{Comparison of total inertia payments [\euro] under different levels of VI installation in the system.}
\label{tab:comp_methods_VI}
\noindent
\centering
    \begin{minipage}{\linewidth} 
    \renewcommand\footnoterule{\vspace*{-5pt}} 
    \begin{center}
        \begin{tabular}{ c || c | c | c }
            \toprule
            \textbf{VI installation} & \textbf{Method 1} & \textbf{Method 2} & \textbf{Method 3} \\ 
            \cline{1-4}
            $0\,\%$ & $258\,160$ & $162\,124$ & $187\,948$ \\
            $10\,\%$ & $245\,195$ & $162\,124$ & $162\,938$ \\
            $20\,\%$ & $156\,725$ & $162\,124$ & $138\,231$ \\
            $35\,\%$ & $135\,849$ & $162\,124$ & $116\,765$ \\
            \arrayrulecolor{black}\bottomrule
        \end{tabular}
        \end{center}
    \end{minipage}
\end{table}

Some interesting observations can be made when studying different levels of virtual inertia in the system provided in Table~\ref{tab:comp_methods_VI}. In particular, as the VI penetration increases the payments from \textit{Method~1} and \textit{Method~3} decrease, while \textit{Method~2} remains unaffected. As a result, the uplift payments become more advantageous for the TSO. Furthermore, by agreeing to pay a fixed amount with the utility function method, this payment scheme can become more expensive as more VI units enter the system. Nonetheless, it is important to keep in mind that the TSO can choose to re-evaluate the inertia payments under such circumstances. Finally, Table~\ref{tab:methods_summary} highlights the main policy implications for all stakeholders involved in the inertia provision. It is clear that all three methods could improve frequency stability and aid the operator in ensuring system reliability. On the other hand, the policy implications to different inertia providers might vary between the payment schemes. Nevertheless, taking all aspects into account, the uplift payments appear to be the most beneficial and practical method of the three. This can be justified on several grounds: (i) the fundamental concepts are already familiar to certain operators; (ii) the method does not require additional \textit{ex-post} calculations; (iii) such payment scheme does not alter the standard UC formulation; and (iv) it prevents excessively high payments to inertia providers. 

\vspace{-0.35cm}
\section{Conclusion}
\label{sec:Conclusion}

This paper proposes three inertia pricing schemes as well as the appropriate methods for reimbursement of respective inertia providers in a transparent and fair manner according to their individual participation. In particular, the focus is on applying the ex-post calculation, utility function and uplift payments, respectively, for determining the appropriate market price for inertial response. A two-step approach based on a frequency-constrained unit commitment formulation is employed, which co-optimizes the provision of energy and inertia services while accounting for their complementary properties and differentiating between the units being online for energy purposes and the ones committed additionally solely for inertia provision. The analysis includes both traditional synchronous generators and converter-based units providing virtual inertia, and gives insights into the impact of each pricing scheme on total system cost and its potential for attracting more affordable providers of inertial response. The results indicate that all three methods could have a beneficial impact on frequency stability and aid the operator in ensuring system reliability. While the policy implications to different inertia providers might vary depending on the payment schemes, the uplift payments appear to be the most beneficial and practical method of the three. This also opens the avenue for future work on integration of inertia provision in the existing ancillary service offering and understanding its impact on other services such as primary frequency response, which could eventually aid the creation of an appropriate market for inertia.

\vspace{-0.35cm}
\bibliographystyle{IEEEtran}
\bibliography{bibliography}

\end{document}